\documentclass{article}

\usepackage{amsthm}
\usepackage{amsmath}
\usepackage{amssymb}

\usepackage{listings}

\lstdefinestyle{Python}{
    language        = Python,
    basicstyle      = \ttfamily,
    keywordstyle    = \color{blue},
    keywordstyle    = [2] \color{teal}, 
    stringstyle     = \color{green},
    commentstyle    = \color{red}\ttfamily
}
\usepackage{graphicx}
\usepackage{framed}
\usepackage{subfigure}
\usepackage{float}

\usepackage{hyperref}

\usepackage{tikz}

\theoremstyle{definition}

\newtheorem{thm}{Theorem}[section]
\newtheorem{lem}[thm]{Lemma}

\newtheorem{defn}[thm]{Definition}
\newtheorem{rem}[thm]{Remark}

\newtheorem*{mainresult}{Main Result}
\newtheorem{heuristic}{Heuristic}

\newcommand{\QED}{\ensuremath{\null\hfill\square}}

\newcommand{\R}{\mathbb{R}}
\newcommand{\Rplus}{\mathbb{R}_{\geq 0}}

\newcommand{\N}{\mathbb{N}}

\newcommand{\dee}{\mathrm{d}}

\newcommand{\ASD}{\text{dim}_\mathcal{F}}

\newcommand{\arctanh}{\mathrm{arctanh}}

\newcommand{\inv}{\mathrm{inv}}

\newcommand{\tsc}{\mathrm{tsc}}
\newcommand{\vol}{\mathrm{vol}}
\newcommand{\diam}{\mathrm{diam}}

\newcommand{\bigO}{\mathrm{O}}

\usepackage{xcolor}

\usepackage{hyperref}
\title{Metric Space Spread, Intrinsic Dimension and the Manifold Hypothesis}

\author{Kevin Dunne\footnote{\url{kevin.dunne[at]mailbox.org}}}

\begin{document}
\maketitle

\begin{abstract}
The concepts of \emph{spread} and \emph{spread dimension} of a metric space were introduced by Willerton in the context of quantifying biodiversity of ecosystems. This paper develops practical applications of spread dimension in the context of machine learning and manifold learning; we show that the topological dimension of a Riemannian manifold can be accurately estimated by computing the spread dimension of a finite subset.

These results are presented as the theoretical basis for a novel method of estimating the intrinsic dimension of data. The practical applications of this method are demonstrated with empirical computations using real and synthetic data.
\end{abstract}

\tableofcontents

\section{Introduction}\label{sec:Introduction}

Willerton introduced the \emph{spread} of a metric space as a measure of ``size'' of a metric space, and a corresponding notion of \emph{spread dimension} based on the growth rate of the spread \cite{Willerton2013:Spread}. Spread is closely related to the \emph{magnitude} of a metric space, which is another measure of ``size'' introduced in the context of quantifying biodiversity of ecosystems \cite{Leinster2013:Magnitude,Leinster2021:Entropy}.

Spread is defined for a broad class of metric spaces, including all finite metric spaces and compact Riemannian manifolds. In this paper we focus on compact Riemannian manifolds and finite subsets of these manifolds, and we show two main results: the spread dimension of a compact Riemannian manifold $M$ coincides with its topological dimension; and, that this value can be approximated from computing the spread dimension of finite subsets $X \subset M$.

We focus on the special case of smooth submanifolds $M \hookrightarrow \R^n$, where the results are presented as the theoretical basis for a novel method of \emph{intrinsic dimension estimation}; a technique to be applied in the context of machine learning and manifold learning.  In Section \ref{sec:Applications} the practical applications of the method are demonstrated with computations using real and synthetic data\footnote{See author's github \url{https://github.com/dk-gh/metric_space_spread} for accompanying source code and data.}.

\subsection{The Spread of a Metric Space}

Every compact metric space $(X,d)$ equipped with a probability measure has a defined spread, which is a strictly positive real value. Since a metric space can be scaled by any constant factor $t \in \Rplus$, spread yields a one-parameter family of values associated with the underlying space.

\begin{defn}
Let $(X,d)$ be a compact metric space equipped with probability measure $\mu$. The \emph{scaled spread} $\sigma_d(t): \Rplus \rightarrow \Rplus$ is the fuction defined by
\begin{equation}\label{eq:definition_of_spread}
\sigma_d(t) = \int_{x\in X} \frac{\dee \mu(x)}{\int_{y\in X}e^{-t d(x,y)}\dee\mu(y)}
\end{equation}
\end{defn}

If $(X,d)$ is finite and equipped with the uniform probability distribution, the scaled spread is given by
\begin{equation}\label{eq:FiniteSpread}
\sigma_{d}(t) = \sum\limits_{x \in X}\frac{1}{\sum\limits_{y \in X}e^{-td(x,y)}}
\end{equation}

For a finite metric space $(X,d)$ the value $\sigma_d(t)$ can be interpreted as encoding the number of distinct points the space resembles when viewed at scale $t$. 
At scales close to zero, the space resembles a single point with $\sigma_d(t)\approx 0$, while at large scales the space resembles $|X|$ discrete points with $\sigma_d(t)\approx |X|$. One of the key insights in \cite{Willerton2013:Spread} is that the growth rate of the spread for intermediate values of $t$ encodes geometric information about the space; it characterises a notion of dimension.

There are a number of ways to define the growth of a function $f(t)$. In \cite{Willerton2013:Spread} Willerton introduced the \emph{instantaneous growth} $\mathcal{G}f(t)$ defined
\begin{equation}\label{eq:definitionGrowthDimension}
\mathcal{G}f(t) = \frac{\dee \ln (f(t))}{\dee \ln(t)} = \frac{t}{f(t)}\frac{\dee}{\dee t}f(t)
\end{equation}
where the second form follows from an application of the chain rule. Willerton defines the spread dimension of a metric sapce in terms of the instantaneous growth of its spread.
\begin{defn}
For $(X,d)$ a metric space with defined spread $\sigma_d(t)$ the \emph{instantaneous spread dimension}, or simply \emph{spread dimension} is defined $\mathcal{G}\sigma_d(t)$.
\end{defn}

Willerton conducted numerical computations on synthetic data which suggest that instantaneous spread dimension characterises a notion of ``dimension at a given scale'' \cite[Section 4]{Willerton2013:Spread}. These computations on finite subsets $X \subset M$ show the value of $\mathcal{G}\sigma_d(t)$ rising towards and plateauing around $\dim(M)$, where $M$ is the space being sampled. This includes examples of fractals where $\dim(M)$ corresponds with the Hausdorff and Minkowski dimension.

Willerton's computations in \cite[Section 4]{Willerton2013:Spread} provide empirical evidence that the spread dimension of a finite subset gives an estimate of the dimension of the space -- this is the core idea developed rigorously in the present work, and in order to make this precise we need another notion of dimension, closely related to the instantaneous spread dimension.

\begin{defn}\label{def:AsymptoticSpreadDimension}
For $(X,d)$ a metric space with defined spread $\sigma_d(t)$, denote the following quantity
\[
\mathcal{F}\sigma_d(t) = \frac{\ln(\sigma_d(t))}{\ln(t)}
\]
and define the \emph{asymptotic spread dimension} $\ASD(X,d)$ as the limit
\begin{equation*}
\ASD(X,d) = \lim\limits_{t \rightarrow \infty}\mathcal{F}\sigma_d(t) 
\end{equation*}
\end{defn}
Note that if the limit $\lim\limits_{t \rightarrow \infty} \mathcal{G}\sigma_d(t)$ exists, then by L'H\^{o}pital's rule the asymptotic spread dimension exists and 
\[\lim\limits_{t \rightarrow \infty} \mathcal{F}\sigma_d(t) = \lim\limits_{t \rightarrow \infty} \mathcal{G}\sigma_d(t)\]
that is, if the asymptotic limit of the instantaneous spread dimension exists, it necessarily coincides with the asymptotic spread dimension.

\begin{rem}
The asymptotic spread dimension is highly reminiscent of a family of fractal dimensions defined in terms of the limit of log-of-size vs log-of-scale. Meckes has shown that the analogous definition for magnitude of a metric space coincides with the Minkowski dimension \cite[Section 7]{Meckes2015:Magnitude}. The instantaneous growth formula (\ref{eq:definitionGrowthDimension})  also appears in a fractal context, namely, in the definition of \emph{Higuchi dimension} \cite{Higuchi1988:AnApproachToAnIrregularTimeSeries}, which is known to approximate the Minkowski dimension \cite{LiehrMassopust2020:OnTheMathematicalValidityOfHiguchi}.
\end{rem}

\subsection{Intrinsic Dimension and the Manifold Hypothesis}

Real-world data is often represented as a finite subset of a high-dimensional space $\R^n$, where each of the $n$-dimensions corresponds with an individual measuring device, for example, a single camera pixel, or an individual question in a questionaire. The high-dimensionality of real-world data is often an artefact of the measurement process, rather than an intrinsic property of the underlying system being measured; points often lie on or close to a lower-dimensional subspace. When data lies on a lower-dimensional subspace with dimension $m$, we say the \emph{intrinsic dimension} of the data is $m$. The intrinsic dimension can also be characterised as the number of parameters needed to explain variance in the data locally. 

There are a number of statistical and machine learning techniques for finding lower-dimensional representations of high-dimensional data. These dimension reduction techniques vary in terms of the underlying formal processes and interpretation, however, in general they can all be characterised as follows: a dimension reduction technique consists of a process taking a subset $X \subset\R^n$ as input, and producing a map $\rho:\R^n \rightarrow \R^k$, for some $k<n$ such that the set $\rho(X) \subset \R^k$ retains some desired statistical or geometric properties of $X$. Dimension reduction techniques can be broadly classified by whether $\rho$ is linear \cite{Cunningham2015:LinearDimensionalityReduction} or non-linear \cite{Maaten2009:DimensionalityReduction}. Non-linear dimension reduction is also known as \emph{manifold learning}. 

A key aspect of manifold learning is the \emph{manifold assumption}, or \emph{manifold hypothesis}: the assumption that data $X \subset \R^n$ lies on or close to a submanifold $M\subset \R^n$. If the manifold hypothesis holds, then the intrinsic dimension of $X$ coincides with the topological dimension of the manifold $M$.

The need to identify the most appropriate target dimension $k$ is a common problem across all forms of dimension reduction. 
The choice of target dimension $k$ for linear dimension reduction methods can be fundamental to an entire fields of study, for example, in psychology \cite{FabrigarEtAl1999:EvaluatingEFA,ZwickVelicer1986:ComparisonofFiveRules}.

In the context of manifold learning, the target dimension $k$ is typically an input parameter for a dimension reduction algorithm. This parameter is then varied and some criterion is chosen to identify the ``correct'' value for $k$. However, there may be more than one reasonable criterion one could choose to optimise for, possibly resulting in different candidate values for $k$.

The importance of identifying the correct target dimension has led to the development of many methods for measuring or estimating the intrinsic dimension of data \cite{Camastra2016:IntrinsicDimension,CampadelliEtAl2015:IntrinsicDimension}.

The results in this paper are presented in the context of justifying spread dimension as a novel technique for intrinsic dimension estimation. These results can be framed in the following terms.

\begin{mainresult}\label{result:1}
For $M\hookrightarrow \R^n$ a submanifold, the topological dimension of $M$ can be estimated from the spread dimension of a finite subset $X \subset M$, where $X$ is a metric space with Euclidean distance inherited from $\R^n$.
\end{mainresult} 

It is crucial for practical applications that this result depends only on the Euclidean distance function, which can be determined from the raw data itself, requiring no additonal knowledge of the geometric properties of $M$.

The theorems proved in Section \ref{sec:DimensionOfManifolds} and Section \ref{sec:FiniteApproximationTheorem} provide the theoretical basis of the Main Result. Two heuristics outlining the use of $\mathcal{G}\sigma_d(t)$ and $\mathcal{F}\sigma_d(t)$ for estimating the intrinsic dimension of data are outlined in Section \ref{sec:Applications}, and we then demonstrate that the method achieves meaningful results in practice when applied to real-world data.

\section{The Spread Dimension of a Manifold}\label{sec:DimensionOfManifolds}

Every Riemannian manifold admits a canonical density, the \emph{Riemannian density} -- see \cite[Proposition 16.45]{Lee2013:SmoothManifolds} for details. This density gives a definition of integration on the manifold, which in turn yields a Radon measure on $M$ -- see Lemma \ref{lem:density_gives_radon} for details. If $M$ is compact then the Radon measure is finite, and hence can be normalised to a probability measure.

Every Riemannian manifold is also a metric space with an intrinsic distance, the \emph{Riemannian distance function} $d_R(x,y)$, defined as the infimum of the lengths of piecewise smooth curves in $M$ joining $x$ and $y$ -- see \cite[Lemma 6.2]{Lee1997:RiemannianManifolds} for details. Therefore every Riemannian manifold has a canonically defined spread $\sigma_{d_R}(t)$.

Reimannian manifolds include the class of smooth submanifolds $M \hookrightarrow \R^m$, in which case $M$ inherits a Riemannian metric from $\R^m$, see, for example \cite[Proposition 13.9]{Lee2013:SmoothManifolds} for details.
For a smooth submanifold $M \hookrightarrow \R^m$ there is another obvious choice of distance function that makes $M$ a metric space -- the usual Euclidean metric $d_E$ inherited from $\R^m$, giving a second natural definition of spread for these spaces, namely $\sigma_{d_E}(t)$, which in general does not coincide with $\sigma_{d_R}(t)$.

In the context of empirical data, one typically has a finite subset $X \subset \R^m$ for some $m$. Under the manifold hypothesis, we assume that $X$ lies on some smooth submanifold $M\subset \R^m$, but we will have no information about the nature of $M$ or its embedding into the Euclidean space, and hence no information about $d_R$.

Some manifold learning techniques attempt to overcome this problem by reconstructing the intrinsic manifold distance function $d_R$ from the distances in Euclidean space, for example, Isomap \cite{TennenbaumEtAl2000:GlobalGeometricFramework}. We take a different approach and prove results using the Euclidean distance $d_E$, rather than reasoning about the Riemannian distance $d_R$. Being able to compute dimension using the Euclidean distance is significant for practical applications, as these distances can be computed directly from empirical data, with no additional assumptions.

\begin{rem}
Non-smooth continuous manifolds are encountered in real data -- for example, piecewise-linear manifolds ocurring in the context of rendering a surface as a polygon mesh -- so it may seem like a strong assumption that data lies on a smooth manifold.
However, by the Stone-Wierstrass Theorem  -- see, for example \cite[Theorem 7.32]{Rudin1976:PrinciplesOfMathematicalAnalysis} -- any continuous embedding of a compact manifold can be uniformly approximated by smooth embeddings. From an empirical point of view, no finite sample could distinguish a smooth from a merely continuous manifold, therefore, there is no practical loss of generality in restricting to smooth manifolds, as opposed to considering manifolds with continuous embeddings $M\hookrightarrow \R^m$.
\end{rem}

In general, it is difficult to compute exact expressions for the spread of a metric space, however, in some cases it is possible to derive an asymptitically equivalent expression. For a pair of functions $f,g: \Rplus \to \Rplus$, we say $f$ and $g$ are \emph{asymptotically equivalent}, denoted $f\sim g$, if 
\[
\frac{f(t)}{g(t)} \to 1 \text{ as } t \to \infty
\]
The following result of Willerton \cite[Theorem 7]{Willerton2013:Spread} gives an asymptotically equivalent expression for the spread $\sigma_{d_R}(t)$ of a Riemannian manifold.
\begin{lem}[Willerton]\label{thm:WillertonManifoldDim}
If is $M$ a compact $n$-dimensional Riemannian manifold where $d_R$ is the Riemannian distance function, then
\begin{equation}\label{eq:manifold_asymp_equiv}
\sigma_{d_R}(t) \sim \frac{1}{n!\omega_n} \bigl(t^n \vol(M) + \frac{n+1}{6}t^{n-2}\tsc(M) + \bigO(t^{n-4})  \bigr)
\end{equation}
where $\vol(M)$ is the volume of $M$, $\tsc(M)$ is the total scalar curvature of $M$, and $\omega_n$ the volume of the unit $n$-ball.
\end{lem}
It will also be useful to consider a weaker notion of asymptotic equivalence, namely, the $\Theta$-class as defined by Knuth \cite{Knuth1976:BigO}.
\begin{defn}
For $f: \Rplus \to \Rplus$ a function, define $\Theta(f)$ as the set of functions $g: \Rplus \to \Rplus$ for which there exist positive constants $A,B$, and $T \in \Rplus$ such that for all $t>T$
\[
Af(t) \leq g(t) \leq Bf(t)
\]

\end{defn}
It is straightforward to verify that $\Theta$ defines an equivalence relation on the set of functions $f: \Rplus \to \Rplus$, and that $f\sim g$ implies $\Theta(f) =\Theta (g)$.

The following lemma shows the asymptotic spread dimension of a metric space can often be inferred from the $\Theta$-class of its scaled spread function.

\begin{lem}\label{lem:RelatingGrowth}
Let $(X,d)$ be a compact metric space.
If $\Theta(\sigma_d(t))=\Theta(t^n)$ for some $n>0$, then $\ASD(X,d)=n$.

\proof{Suppose $\Theta(\sigma_d(t))=\Theta(t^n)$ then there exists $A,B >0$ and $T \in \Rplus$ such that for all $t> T$
\[
At^n \leq \sigma_d(t) \leq Bt^n
\]

We can choose $T$ such that $T>1$, and therefore, for all $t> T$ we have
\[
\frac{\ln(At^n)}{\ln(t)} \leq \frac{\ln(\sigma_d(t))}{\ln(t)}  \leq \frac{\ln(Bt^n)}{\ln(t)} 
\]

Since $\ln(At^n) = n\ln(t)+\ln(A)$, then for all $t>T$ we have
\[
n + \frac{\ln(A)}{\ln(t)}\leq \frac{\ln(\sigma_d(t))}{\ln(t)}  \leq n + \frac{\ln(B)}{\ln(t)}
\]
and hence
\[
 \frac{\ln(\sigma_d(t))}{\ln(t)} \rightarrow n \text{ as } t \rightarrow \infty
\]
as required.
\QED}
\end{lem}

From the previous lemma we can conclude that the asymptotic spread dimension of an $n$-dimensional Riemannian manifold is $n$.

\begin{thm}\label{thm:IntrinsicMetricSpreadDimension}
Let $M$ be an $n$-dimensional Riemannian manifold with $d_R$ the Riemannian distance function, then $\ASD(M,d_R) = n$.
\proof{The expression on the right hand side of (\ref{eq:manifold_asymp_equiv}) is asymptotically equivalent to $\frac{\vol(M)}{n!\omega_n} t^n$. Since $\sim$ is an equivalence relation, Lemma \ref{thm:WillertonManifoldDim} implies
\[
\sigma_{d_R}(t) \sim \frac{\vol(M)}{n!\omega_n} t^n
\]
and hence we have
\[\Theta(\sigma_{d_R}(t)) = \Theta\Big(\frac{\vol(M)}{n!\omega_n} t^n\Big) = \Theta(t^n)\]
and therefore the required result follows from Lemma \ref{lem:RelatingGrowth}. \QED}
\end{thm}

\subsection{The Spread Dimension of Euclidean Submanifolds}

Extending Theorem \ref{thm:IntrinsicMetricSpreadDimension} to the case of the Euclidean distance function will require that we show $d_E$ and $d_R$ are bilipschitz, and that asymptotic spread dimension is invariant under bilipschitz maps.

\begin{defn}
For $(X,d_X)$ and $(Y,d_Y)$ metric spaces, a function $f: X \to Y$ is \emph{Lipschitz continuous} if there exists $A \in \Rplus$ such that for all, $x,x'\in X$ the following inequality holds
\[
d_Y\big(f(x),f(x') \big)\leq Ad(x,x')
\]
A function $f: X \to Y$ is \emph{bilipschitz} if there exists $A,B \in \Rplus$ such that for all $x,x' \in X$
\[
A d_X(x,x') \leq d_Y(f(x), f(x')) \leq B d_X(x,x') 
\]

For two metrics $d_1$ and $d_2$ defined on the set $X$, we say that $d_1$ and $d_2$ are \emph{bilipschitz} if the identity map $\text{id} : (X,d_1) \rightarrow (X,d_2)$ is bilipschitz.
\end{defn}

The following result states that asymptotic spread dimension is preserved by bilipschitz maps.

\begin{lem}\label{lem:BilipschitzLemma}
Let $(X,d_X)$ and $(Y,d_Y)$ be compact metric spaces such that there exists a bilipschitz map $f:X\to Y$. If $\ASD(X,d_X)$ exists then so does $\ASD(Y,d_Y)$ and $\ASD(X,d_X)=\ASD(Y,d_Y)$.

\proof{ Suppose $\ASD(X,d_X)= n$, and that $d_X$ and $d_Y$ are bilipschitz, i.e. there exists $A,B >0 $ such that for all $x,y \in X$
\begin{equation}\label{eq:bilipleminequalities}
Ad_X(x,y) \leq d_Y(x,y) \leq B d_X(x,y) 
\end{equation}

From the definition of $\sigma_d(t)$ it is easy to see that for any $C >0 $ we have $\sigma_{d_X}(Ct) = \sigma_{Cd_X}(t)$. Using the fact that Lebesgue integration is order preserving \cite[Theorem 3.22]{Axler2020:MeasureTheory}, and by unpacking the definition of $\sigma_d(t)$, the inequalities (\ref{eq:bilipleminequalities}) imply
\begin{equation}\label{eq:bilipschitzInequalities}
\sigma_{d_X}(At) \leq \sigma_{d_Y}(t) \leq \sigma_{d_X}(Bt) 
\end{equation}
Consider the first inequality of (\ref{eq:bilipschitzInequalities}). It follows immediately that for all $t >1$
\[
\frac{\ln(\sigma_{d_X}(At))}{\ln(At)} \leq \frac{\ln(\sigma_d(t))}{\ln(At)} 
\]
By assumption $\ASD(X,d_X)= n$, hence we have $\lim\limits_{t \rightarrow \infty} \frac{\ln(\sigma_{d_X}(At))}{\ln(At)} = n$.

Note that
\begin{align*}
\lim\limits_{t \rightarrow \infty }\frac{\ln(\sigma_{d_Y}(t))}{\ln(At)} &= \lim\limits_{t \rightarrow \infty }\frac{\ln(\sigma_{d_Y}(t))}{\ln(t) + \ln(A)} \\
 &= \lim\limits_{t \rightarrow \infty }\frac{\ln(\sigma_{d_Y}(t))}{\ln(t)}
\end{align*}
Hence as $t \rightarrow \infty$ we have
\begin{equation}\label{eq:proof_one_way}
n \leq \frac{\ln(\sigma_{d_Y}(t))}{\ln(t)}
\end{equation}
By the same argument using the second inequality from (\ref{eq:bilipschitzInequalities}) we obtain the reverse inequality to (\ref{eq:proof_one_way}), and hence we conclude that $\lim\limits_{t \rightarrow \infty}\frac{\ln(\sigma_{d_Y}(t))}{\ln(t)}$ exists and is equal to $n$, as required.
 \QED}
\end{lem}

\begin{rem}\label{rem:BilipschitzFractalDim}
Analogous to Lemma \ref{lem:BilipschitzLemma} which shows asymptotic spread dimension is preserved by bilipschitz maps, the \emph{Hausdorff dimension} and \emph{packing dimension} of a metric space is also preserved by bilipschitz maps, see, for example \cite[Theorem 6.8.9]{Edgar2008:MeasureTopologyAndFractalGeometry}.
\end{rem}

To show that the Euclidean and Riemannian distance functions are bilipschitz, we rely on some results concerning \emph{tubular neighbourhoods} of smooth submanifolds of Euclidean space. For $i:M \hookrightarrow \R^m$ an $n$-dimensional submanifold, the tangent space $T_xM$ at each $x\in M$ can be canonically identified with a subspace of $T_x\R^m$. At each point $x \in M$ we can therefore define the \emph{normal space} $N_xM$ to be the $(m-n)$-dimensional space consisting of those vectors in $T_x\R^m$ but not in $T_xM$. The \emph{normal bundle} $NM$ is defined to be the union $\bigsqcup\limits_{x \in M} N_xM$ with points of the form $(x,v)$ where $x \in M$ and $v \in N_xM$. For more details see \cite[pp. 137-138]{Lee2013:SmoothManifolds}.

\begin{defn}
Let $i:M \hookrightarrow \R^m$ be a submanifold. A \emph{tubular neighbourhood} of $M$ is an open neighbourhood $U\subset\R^m$ containing $i(M)$, such that there exists a subspace $V\subset NM$ of the form
\[
V = \big\{\, (x,v) \, \big| \, x \in M, \, |v|< c(x), \text{ where } c \text{ is continuous } \big\}
\]
where $U$ is diffeomorphic to $V$.
\end{defn}
The following can be found in \cite[Theorem 6.24]{Lee2013:SmoothManifolds}.
\begin{lem}\label{lem:TubularNeighbourhoodsExist}
Every smooth submanifold $i:M \hookrightarrow \R^m$ has a tubular neighbourhood.
\end{lem}
The following can be found in \cite[Theorem 6.25]{Lee2013:SmoothManifolds}.
\begin{lem}\label{lem:SmoothRetractExists}
Let $i:M \hookrightarrow \R^m$ be a smooth submanifold. For any tubular neighbourhood $U$ of $M$, there exists a smooth map $\rho: U \rightarrow M$ that $\rho(i (x)) = x$ for all $x \in M$.
\end{lem}

The smooth map $\rho: U \rightarrow M$ is obtained by composing the diffeomorphism $\tau: U \rightarrow V$ with the projection map $\pi:NM \to M$, defined by $(x,v) \mapsto x$.

\begin{lem}\label{lem:ManifoldMetricsAreBilipschitz}
Let $i:M \hookrightarrow \R^m$ be a smooth submanifold. The distance functions $d_R$ and $d_E$ are bilipschitz, where $d_E$ is the Euclidean distance function restricted to $M$, and let $d_R$ is the Riemannian distance function on $M$.
\proof{It is clear that for all $x,y \in M$ the inequality $d_E(i(x),i(y)) \leq d_R(x,y)$ holds, and hence all that remains to  show is there exists $K \in \Rplus$ such that $d_R(x,y) \leq Kd_E(x,y)$ for all $x,y\in M$.

By Lemma \ref{lem:TubularNeighbourhoodsExist} there exists a tubular neighbourhood $U \subset \R^m$ which inherits the Euclidean distance function $d_E$ from $\R^m$.
By Lemma \ref{lem:SmoothRetractExists} there exists a smooth map $\rho: U \rightarrow M$ satisfying $\rho(i(x)) = x$ for all $x \in M$.

Since $\rho$ is smooth it is Lipschitz continuous -- see, for example \cite[Proposition C.29]{Lee2013:SmoothManifolds} -- hence, there exists $K\in\Rplus$ such that 
\[
d_R(\rho(u),\rho(v)) \leq K d_E(u, v)
\]
and in particular this holds for those points in $U$ of the form $u = i(x)$, and $v = i(y)$, in which case we have
\begin{align*}
d_R(x,y) &= d_R(\rho(i(x)),\rho(i(y))) \\
&\leq K d_E(i(x), i(y))
\end{align*}
as required. \QED}
\end{lem}

We can now show the Euclidean version of Theorem \ref{thm:IntrinsicMetricSpreadDimension}.

\begin{thm}\label{thm:ASDofManifold}
If $(M,d)$ is an $n$-dimensional manifold with a smooth embedding $M \hookrightarrow \R^m$ for some $m$, where $d_E$ is the Euclidean distance function of $\R^m$ restricted to $M$, then $\ASD(M,d_E) = n$.
\proof{By Theorem \ref{thm:IntrinsicMetricSpreadDimension} we have $\ASD(M,d_R) = n$, where $d_R$ is the Riemannian distance function on $M$. By Lemma \ref{lem:ManifoldMetricsAreBilipschitz} the distance functions $d_E$ and $d_R$ are bilipschitz, and hence by Lemma \ref{lem:BilipschitzLemma} we have $\ASD(M,d_E) = n$, as required. \QED}
\end{thm}

\subsection{The Local Neighbourhood Dimension}

We will now show a local version of Theorem \ref{thm:ASDofManifold} which states that around each point in an $n$-dimensional manifold, there exists a local neighbourhood with asymptotic spread dimension $n$. This result essentially follows from the fact that Euclidean $n$-balls have asymptotic spread dimension $n$, which is demonstrated below. From an empirical point of view, this result is useful as it justifies probing the dimension of dataset by considering local samples around specific points.

Willerton gave the following characterisation of the spread of the unit interval, and evaluated its asymptotic behaviour \cite[Theorem 5]{Willerton2013:Spread}.

\begin{lem}[Willerton]\label{lem:WillertonSpreadOfInterval}
If $I = [0,1]$ is the unit interval equipped with the usual distance function $d(x,y) = |x-y|$ then
\[
\sigma_d(t) = \frac{\arctanh\bigr(\sqrt{1 - e^{-t}}\bigl)}{\sqrt{1 - e^{-t}}}
\]
and the asymptotic behaviour of this expression is characterised by
\begin{equation}\label{eq:asymptoticInterval}
\Big|  \sigma_d(t) - \Big(\frac{t}{2} + \ln(2) \Big)  \Big| \rightarrow 0 \text{ as } t \rightarrow \infty
\end{equation}
\end{lem}

\begin{lem}\label{lem:SpreadDimensionInterval}
For $I=[0,1]$ the unit interval equipped with the usual distance function $d$, we have $\ASD(I,d) = 1$.
\proof{Using the formula for $\sigma_d(t)$ from Lemma \ref{lem:WillertonSpreadOfInterval} one can compute the limit $\ASD(I,d)$ directly. 

Alternatively, it is easy verify from (\ref{eq:asymptoticInterval}) that $\Theta(\sigma_d(t)) = \Theta(t)$, hence the result follows immediately from Lemma \ref{lem:RelatingGrowth}
\QED}
\end{lem}

For $(X,d_X)$ and $(Y,d_Y)$ metric spaces, let $(X\times Y, d_X \oplus d_Y)$ denote the product space, where $d_X \oplus d_Y$ is defined
\[
d_X \oplus d_Y((x_1,y_1),(x_2,y_2)) = d_X(x_1,x_2) + d_Y(y_1, y_2)
\]

\begin{lem}\label{lem:SpreadProduct}
If $(X,d_X)$ and $(Y,d_Y)$ are compact metric spaces equipped with probability measures, then
\[
\ASD(X\times Y, d_X \oplus d_Y)  = \ASD(X,d_X) + \ASD(Y,d_Y)
\]

and
\[
\mathcal{G}\big(\sigma_{d_X \oplus d_Y}(t)\big) = \mathcal{G}\big(\sigma_{d_X}(t)\big)  + \mathcal{G}\big(\sigma_{d_Y}(t)\big) 
\]
\proof{The results follow from Fubini's Theorem -- see, for example \cite[Theorem 5.32]{Axler2020:MeasureTheory}, which is applied twice, as follows

\begin{align*}
\sigma_{d_X \oplus d_Y}(t) &= \int\limits_{(x_2,y_2)\in X\times Y} \frac{\dee\mu(x_2,y_2)}{\int\limits_{(x_1,y_1)\in X\times Y} e^{-t d_X \oplus d_Y((x_1,y_1),(x_2,y_2))} \dee\mu(x_1,y_1)} \\
&= \int\limits_{(x_2,y_2)\in X\times Y} \frac{\dee \mu(x_2)\dee\mu(y_2)}{\Bigr( \int\limits_{x_1\in X}  e^{-t d_X(x_1,x_2)} \dee\mu(x_1)\Bigl) \Bigr(\int\limits_{y_1\in Y} e^{-t  d_Y(y_1, y_2)} \dee \mu(y_1)\Bigl) } \\
&= \int\limits_{x_2\in X} \frac{\dee \mu(x_2)}{\int\limits_{x_1\in X}  e^{-t d_X(x_1,x_2)} \dee \mu(y_1) } \int\limits_{y_2\in Y} \frac{\dee\mu(y_2)}{\int\limits_{y_1\in Y} e^{-t d_Y(y_1, y_2)} \dee \mu(y_1)}\\
&= \sigma_{d_X}(t) \sigma_{d_Y}(t)
\end{align*}

The first result then follows immediately from
\[
\frac{\ln\big( \sigma_{d_X}(t)  \sigma_{d_Y}(t) \big)}{\ln(t)} = \frac{\ln\big( \sigma_{d_X}(t) \big)}{\ln(t)}  + \frac{\ln\big(\sigma_{d_Y}(t) \big)}{\ln(t)} 
\]
The second result follows from $\mathcal{G}\big(\sigma_{d_X}(t)  \sigma_{d_Y}(t)\big) = \mathcal{G}\big(\sigma_{d_X}(t) \big)+ \mathcal{G}\big(  \sigma_{d_Y}(t)\big)$, which is a straightforward application of the product rule.
\QED}
\end{lem}

It follows immediately from Lemma \ref{lem:SpreadProduct} that if $\ASD(X,d_X) = n$ and $\ASD(Y,d_Y) = m$ then $\ASD(X\times Y, d_X \oplus d_Y) = n+m$.

\begin{rem}\label{rem:FractalProducts}
While $\dim(X \times Y) = \dim(X) + \dim(Y)$ holds for topological dimension of Euclidean spaces and manifolds, it does not hold for topological spaces in general, or other definitions of dimension. In particular, Hausdorff and packing dimension do not satisfy this in general \cite[Lemma 2.2]{Xiao1996:CartesianProduct}. 
\end{rem}

\begin{lem}\label{lem:BallASDisn}
For each $n\geq 1$ and $r>0$ let $B_r(0) \subset \R^n$ be the $n$-ball of radius $r$. If $d_E$ is the usual Euclidean distance, then $\ASD(B_r(0),d_E)=n$.
\proof{Note that for any $r \in\Rplus$ the ball $B_r(0)$ is bilipschitz equivalent to the unit ball $B_1(0)$ under scaling by $r$, hence by Lemma \ref{lem:BilipschitzLemma} it is enough to show that $\ASD(B_1(0),d_E)=n$.

By Lemma \ref{lem:SpreadDimensionInterval} we have $\ASD(I,d) = 1$, hence by Lemma \ref{lem:SpreadProduct} we have $\ASD(I^n,d_p) =n$ where $d_p$ is the product metric.

It is straightforward to show that the product metric distance function $d_p$ and the Euclidean distance function $d_E$ are bilipschitz, and therefore by Lemma \ref{lem:BilipschitzLemma} we have $\ASD(I^n,d_E) =n$.

For each $n\in \N$, the unit $n$-cube and unit $n$-ball $B_1(0)$ equipped with the usual Euclidean distance functions are also bilipschitz -- see, for example \cite[Corollary 3]{GriepentrogEtAl2008:ABiLipschitzContinuous}, hence $\ASD(B_1(0),d_E)=n$, as required. 
\QED}
\end{lem}

To prove the local dimension theorem we will use the fact that each point in $M$ has a \emph{normal neighbourhood}, see, for example \cite[Chapter 5]{Lee1997:RiemannianManifolds}. Recall a normal neighbourhood of a point $x \in M$ is a neighbourhood $U$ of $x$ such that $U$ is diffeomorphic to a subspace of the tangent space $V \subset T_x M$ under the exponential map.

If the subspace $V \subset T_xM$ is of the form $V = B_r(0)$ for some $r>0$ then the neighbourhood $U$ is called a \emph{geodesic ball}. Every geodesic ball $U\subset M$ around a point $x \in M$ is of the form $U = \{ y \, | \, d_R(x,y) < r \}$, see for example \cite[Corollary 6.11]{Lee1997:RiemannianManifolds}.

\begin{thm}\label{thm:DimensionOfNeighbourhood}
Let $M$ be an $n$-dimensional smooth manifold. For every $x\in M$ there exists a neighbourhood $U\subset M$ such that the asymptotic spread dimension of $U$ is $n$.

\proof{ Around each point $x \in M$ there exists a geodesic ball $U$. The geodesic ball $U$ is diffeomorphic to an $n$-ball $B_r(0)$ for some $r>0$, and therefore bilipschitz equivalent to $(B_r(0), d_E)$. Lemma \ref{lem:BallASDisn} states that $\ASD(B_r(0),d_E) = n$ and hence by Lemma \ref{lem:BilipschitzLemma} we also have $\ASD(U,d_R) = n$, as required.
\QED}
\end{thm}

\section{Estimating Spread Dimension}\label{sec:FiniteApproximationTheorem}

The main result of this section is that the asymptotic spread dimension of a Riemannian manifold can be approximated from its finite subsets.

The following definition is a condition whereby a space can be divided into arbitrarly many equally sized parts, each of which can be made arbitrarily small in the sense of diameter $\diam(U) = \sup\{d(x,y) | x,y \in U \}$.

\begin{defn}\label{def:FiniteSamplingCondition}
Let $(X, d)$ be a compact metric space equipped with a probability measure $\mu$, and let $P_n = \{ I^n_i\}$ be a sequence of partitions of $X$. We say that $X$ and $P_n$ satisfy the \emph{finite sampling condition} if the following conditions hold
\begin{enumerate}
\item for each $n\in \N$ the partition satisfies $|P_n| = n$, and $\mu(I^n_i) = \frac{1}{n}$ for every component $I^n_i \in P_n$;
\item and, for all $\varepsilon>0$ there exists $N \in \N$ such that for all $n>N$, $\diam(I^n_i)<\varepsilon$ for all $I^n_i\in P_n$.
\end{enumerate}
The space $(X,d)$ is said to satisfy the finite sampling condition if such a sequence $P_n$ exists. A set of \emph{representative points} $X_n$ for the partition is a set with $|X_n| = n$ where each $x \in X_n$ belongs to exactly one element of the partition $P_n$.
\end{defn}

The intuition behind this definition is that the sequence of partitions acts like a uniform mesh covering the space which can be made arbitrarily fine. A random sample therefore will resemble a set of representative points with respect to this uniform mesh.

\subsection{The Spread Dimension Estimation Theorems}

Recall, for $P = \{A_1,...,A_n \} $ a partition of $X$ and $f: X \to \Rplus$ a measurable function, the \emph{lower Lebesgue sum} $\mathcal{L}(f,P)$ and \emph{upper Lebesgue sum} $\mathcal{U}(f,P)$ are defined

\[
\mathcal{L}(f,P) = \sum\limits_{i=1}^n \mu(A_i) \inf\limits_{A_i} f, \quad\qquad
\mathcal{U}(f,P) = \sum\limits_{i=1}^n \mu(A_i) \sup\limits_{A_i} f
\]

The Lebesgue integral is typically defined as the supremum of lower Lebesgue sums over all measurable partitions

\[
\int_{X} f \dee\mu = \sup \{ \, \mathcal{L}(f,P) \,\, | \,\, P \text{ a partition of } X \, \}
\]

However, if $\mu(X)<\infty$ and $f$ is bounded then the Lebesgue integral coincides with the infimum of upper Lebesgue sums \cite[p. 99]{Axler2020:MeasureTheory}
\[
\int_{X} f \dee\mu = \inf \{ \, \mathcal{U}(f,P) \,\, | \,\, P \text{ a partition of } X \, \}
\]

\begin{lem}\label{lem:ApproximationLemma}
Let $X$ be a compact metric space and $P_n={I_i^n}$ a sequence of partitions satisfying the finite sampling condition, with $X_n \subset X$ a set of representative points. If $f:X \rightarrow \R$ is a continuous function with $f_n: X_n \rightarrow \R$ defined by the restriction $f_n = f|_{X_n}$, then
\[
\sum\limits_{x\in X_n} f_n(x)\frac{1}{n} \rightarrow \int\limits_{x\in M} f(x) \dee\mu(x) \quad\text{as}\quad n \rightarrow \infty
\]

\proof{For all $n \in \N$ we have both
\[
\mathcal{L}(f, P_n) \leq \sum\limits_{x\in X_n} f_n(x)\frac{1}{n} \leq \mathcal{U}(f, P_n)
\]
and
\[
\mathcal{L}(f, P_n) \leq \int\limits_{x\in M} f(x) \dee\mu(x)\leq \mathcal{U}(f, P_n)
\]
hence it is enough to show that
\begin{equation*}
| \mathcal{L}(f, P_n) - \mathcal{U}(f, P_n)| \rightarrow 0 \quad\text{as}\quad n \rightarrow \infty
\end{equation*}

Note that since $X$ is compact and $f$ is continuous, by the Heine-Cantor Theorem -- see, for example \cite[Proposition 5.8.2]{Sutherland1975:IntroductionMetricTopology} -- $f$ is uniformly continuous, and hence for any $\varepsilon>0$, there exists $\delta > 0$ such that for all $x,y \in X$
\begin{equation}\label{eq_little_eq}
d(x,y)< \delta \quad \Rightarrow \quad |f(x) - f(y)| <\varepsilon
\end{equation}

By the finite sampling condition we can find $N \in \N$ such that for all $n >N$ we have $\diam(A^n_i)<\delta$, for all $A^n_i \in P_n$, and hence for all $n>N$ we have
\begin{align*}
| \mathcal{L}(f, P_n) - \mathcal{U}(f, P_n)| &= \Big| \sum\limits_{i=1}^n \frac{1}{n} \inf\limits_{A^n_i} f - \sum\limits_{i=1}^n \frac{1}{n} \sup\limits_{A^n_i} f \Big| \\
&=\Big| \sum\limits_{i=1}^n \frac{1}{n} (\inf\limits_{A^n_i} f -  \sup\limits_{A^n_i} f) \Big| \\
&\leq \sum\limits_{i=1}^n \frac{1}{n} \Big|(\inf\limits_{A^n_i} f -  \sup\limits_{A^n_i} f)\Big|  \\
&< \sum\limits_{i=1}^n \frac{1}{n} \varepsilon  \quad\text{ by (\ref{eq_little_eq})} \\
& = \varepsilon
\end{align*}
as required. \QED}
\end{lem}

The following result can be found in, for example, \cite[Theorem 9.16]{Apostol1974:Analysis}.
\begin{lem}\label{lem:ApostolLemma}
Let $f: \N \times \N \rightarrow \R$ be a double sequence, and let $g_n(m)= f(n,m)$. If the convergence $g_n \rightarrow g$ as $n \rightarrow \infty$ is uniform, and the limit $\lim\limits_{m \rightarrow \infty}g(m)=L$ exists, then the double limit $\lim\limits_{(n,m)\rightarrow \infty} f(n,m)$ exists and is equal to $L$. 
\end{lem}

The following lemma shows that the spread of a metric space can be approximated from the spread of finite subsets.

\begin{lem}\label{lem:SpreadApproximation}
Let $(M,d)$ be a compact metric space with a sequence of partitions $P_n$ satisfying the finite sampling condition. For each $n\in \N$ let $X_n \subset M$ be a set of representative points consisting of one point $x_i \in I_i^n$ for each subset in the partition $I_i^n \in P_n$. For each $t \in\Rplus$
\[
\sigma_{d_n}(t) \rightarrow \sigma_d(t) \text{ as } n \rightarrow \infty
\]
where $\sigma_{d_n}(t)$ denotes the spread of $X_n$.

\proof{
For each $x \in M$ let 
\[\varphi_n(x) =\sum\limits_{y \in X_n} e^{-td(x,y)} \frac{1}{n} \quad\text{ and }\quad\varphi(x) = \int\limits_{y \in M} e^{-td(x,y)} \dee\mu(y)\]
By Lemma \ref{lem:ApproximationLemma} we have $\varphi_n(x) \to \varphi(x)$ as $n \to \infty$. We will now show that this convergence is uniform in $x$.

For each fixed $t \in\Rplus$ the function $e^{-td(x,y)}:M \times M\rightarrow \Rplus$ is continuous. Since $M$ is compact, then by the Heine-Cantor Theorem -- see, for example \cite[Proposition 5.8.2]{Sutherland1975:IntroductionMetricTopology} -- the function is uniformly continuous, that is, for all $\varepsilon>0$, there exists $\delta>0$ such that for all $x,y_0,y_1 \in M$ we have
\begin{equation}\label{eq_intermediate_result}
d(y_0,y_1)<\delta \quad\Rightarrow \quad \big| e^{-td(x,y_0)} -  e^{-td(x,y_1)}   \big| < \varepsilon
\end{equation}

Hence for all $\varepsilon>0$, there exists $N\in \N$ such that for all $n>N$, we have $\diam(A^n_i)<\delta$ such for all $A^n_i \in P_n$, and therefore, for all $x\in M$ we have

\begin{align*}
\big| \varphi_n(x) - \varphi(x) \big| &\leq | \mathcal{L}(e^{-td(x,y)}, P_n) - \mathcal{U}(e^{-td(x,y)}, P_n)| \\
&\leq \sum\limits_{i=1}^n \frac{1}{n} \Big|\big(\inf\limits_{y_0 \in A^n_i} e^{-td(x,y_0)} -  \sup\limits_{y_1 \in A^n_i} e^{-td(x,y_1)}\big)\Big|  \\
&\leq \sum\limits_{i=1}^n \frac{1}{n} \varepsilon \qquad \text{by (\ref{eq_intermediate_result})} \\
&= \varepsilon
\end{align*}

Hence the convergence $\varphi_n(x) \rightarrow \varphi(x)$ as $n \to \infty$ is uniform in $x$. We will now show that the following also converges uniformly in $x$
\[
\frac{1}{\varphi_n(x)} \to \frac{1}{\varphi(x)} \quad\text{ as }\quad n \to \infty
\]

Note that $\varphi: M \to \R$. Let $r =\diam(M)$. Hence for all $x,y \in M$ we have:

\begin{align*}
\varphi(x) &= \int\limits_{y \in M} e^{-td(x,y)} \dee\mu(y) \geq \int\limits_{y \in M} e^{-tr} \dee\mu(y) =e^{-tr}\mu(M) =e^{-tr}
\end{align*}
and similarly, for each $n \in \N$ we have
\begin{align*}
\varphi_n(x) &= \sum\limits_{y \in X_n} e^{-td(x,y)} \frac{1}{n} \geq \sum\limits_{y \in X_n} e^{-tr} \frac{1}{n} =ne^{-tr}\frac{1}{n}=e^{-tr} \\
\end{align*}
Denoting $L = e^{-t r}$ we have
\[
\Big|\frac{1}{\varphi_n(x)} - \frac{1}{\varphi(x)} \Big| = \Big| \frac{\varphi(x)-\varphi_n(x)}{\varphi_n(x)\varphi(x)} \Big|
 \leq \Big| \frac{\varphi(x)-\varphi_n(x)}{L^2} \Big| 
\]
Since $\varphi_n(x) \rightarrow \varphi(x)$ uniformly in $x$ as $n\to \infty$, for each $\varepsilon > 0$ we can pick $N\in \N$ such that for all $m > N$ we have $|\varphi_m(x) - \varphi(x)| < \varepsilon L^2$. Therefore for all $\varepsilon>0$ there exists $N \in \N$ such that for all $m> \N$
\[
\Big|\frac{1}{\varphi_n(x)} - \frac{1}{\varphi(x)} \Big|< \varepsilon
\]
and therefore we have uniform convergence $\frac{1}{\varphi_n(x)} \to \frac{1}{\varphi(x)}$ as $n \to \infty$.

Consider the double sequence
\[
f(n,m) = \sum\limits_{x \in X_m} \Big(\frac{1}{\varphi_n(x)} \Big)\frac{1}{m}
\] 
and let $g_n(m) = f(n,m)$, and let $g(m) = \lim\limits_{n \to \infty} g_n(m)$

First we show that the iterated limit is equal to the spread of $\sigma_d(t)$

\begin{align*}
\lim\limits_{m \rightarrow \infty} \lim\limits_{n \rightarrow \infty}f(m,n) &= \lim\limits_{m \rightarrow \infty} \lim\limits_{n \rightarrow \infty} \sum\limits_{x \in X_m} \Big(\frac{1}{\varphi_n(x)} \Big)\frac{1}{m} \\
&= \lim\limits_{m \rightarrow \infty} \sum\limits_{x \in X_m} \lim\limits_{n \rightarrow \infty} \Big(\frac{1}{\varphi_n(x)} \Big)\frac{1}{m} \\
&= \lim\limits_{m \rightarrow \infty} \sum\limits_{x \in X_m} \Big( \frac{1}{\varphi(x)}\Big)\frac{1}{m} \\
&= \int\limits_{x \in M} \frac{\dee\mu(x)}{\varphi(x)} \\
&= \int\limits_{x \in M} \frac{\dee\mu(x)}{\int\limits_{y \in M} e^{-td(x,y)}\dee\mu(y)}
\end{align*}
In order to apply Lemma \ref{lem:ApostolLemma}, we will now show that the convergence $g_n(m) \to g(m)$ as $n \to \infty$ is uniform in $m$. Note that by the uniform convergence of $\frac{1}{\varphi_n(x)} \to \frac{1}{\varphi(x)}$ as $n \to \infty$, for all $\varepsilon> 0$ there exists $N \in \N$ such that for all $n >N$ and for all $x \in M$ we have

\begin{equation}\label{eq:UseUniformConvergence}
\Big|\frac{1}{\varphi_n(x)} - \frac{1}{\varphi(x)} \Big| < \varepsilon
\end{equation}

Now consider
\begin{align*}
|g_n(m) - g(m)| &= \Big|\Big( \sum\limits_{x \in X_m} \frac{1}{\sum\limits_{y \in X_n} e^{-td(x,y)}\frac{1}{n}}\frac{1}{m}\Big) - \Big( \sum\limits_{x \in X_m} \frac{1}{\int\limits_{y \in M} e^{-td(x,y)}\dee\mu(y)}\frac{1}{m}\Big) \Big|\\
&=\Big| \sum\limits_{x \in X_m}\Big( \frac{1}{\sum\limits_{y \in X_n} e^{-td(x,y)}\frac{1}{n}}- \frac{1}{\int\limits_{y \in M} e^{-td(x,y)}\dee\mu(y)}\Big)\frac{1}{m} \Big| \\
&\leq \sum\limits_{x \in X_m}\Big| \Big( \frac{1}{\sum\limits_{y \in X_n} e^{-td(x,y)}\frac{1}{n}}- \frac{1}{\int\limits_{y \in M} e^{-td(x,y)}\dee\mu(y)}\Big)\frac{1}{m} \Big| \\
&< \sum\limits_{x \in X_m} \frac{\varepsilon}{m} \qquad\text{by (\ref{eq:UseUniformConvergence})} \\
&= \varepsilon
\end{align*}
and therefore the convergence $g_n(m) \to g(m)$ as $n\to \infty$ is uniform in $m$, and we can apply Lemma \ref{lem:ApostolLemma} to conclude that the iterated limit above is equal to the double limit
\[\lim\limits_{(n,m) \rightarrow \infty}f(m,n) = \int\limits_{x \in M} \frac{\dee\mu(x)}{\int\limits_{y \in M} e^{-td(x,y)}\dee\mu(y)}\]
In particular we have the limit when $n=m$ 
\[
\sum\limits_{x \in X_n} \frac{1}{\sum\limits_{y \in X_n} e^{-td(x,y)}\frac{1}{n}}\frac{1}{n} = \sum\limits_{x \in X_n} \frac{1}{\sum\limits_{y \in X_n} e^{-td(x,y)}}\]
and therefore, for each $t \in \Rplus$ we have
\begin{align*}
\lim\limits_{n \to \infty}\sigma_{d_n}(t) &= \lim\limits_{n \rightarrow \infty}\sum\limits_{x \in X_n} \frac{1}{\sum\limits_{y \in X_n} e^{-td(x,y)}} 
= \int\limits_{x \in M} \frac{\dee\mu(x)}{\int\limits_{y \in M} e^{-td(x,y)}\dee\mu(y)} 
\end{align*}
as required.
\QED}
\end{lem}

\begin{rem}
Generally, the convergence $\sigma_{d_n}(t) \rightarrow \sigma_d(t)$ is not uniform in $t$. That is, for any fixed $t$, the value $\sigma_d(t)$ can be approximated arbitrarily closely by finite subsets, however, increasing the value of $t$ may increase the number of points required to obtain a close approximation. The practical significance of this is illustrated in Section \ref{sec:Applications}.
\end{rem}

It follows more of less immediately that the quantity $\mathcal{F}\sigma_d(t)$ from Definition \ref{def:AsymptoticSpreadDimension} can also be approximated from finite subsets.

\begin{thm}\label{thm:FiniteApproximation1}
Let $(M,d)$ be a compact metric space with a sequence of partitions $P_n$ satisfying the finite sampling condition. For each $n\in \N$ let $X_n \subset M$ be a set of representative points consisting of one point $x^n_i \in I_i^n$ for each subset in the partition $I_i^n \in P_n$. For each $t \in\Rplus$
\[
\mathcal{F}\sigma_{d_n}(t)\rightarrow \mathcal{F}\sigma_d(t)\quad\text{ as }\quad n \rightarrow \infty
\] 
where $\sigma_{d_n}(t)$ denotes the spread of $X_n$.

\proof{By Lemma \ref{lem:SpreadApproximation} we have $\sigma_{d_n}(t)\rightarrow \sigma_d(t)$ as $ n \rightarrow \infty$, and since $\ln(x)$ is continuous we have
\[
\frac{\ln(\sigma_{d_n}(t))}{\ln(t)} \rightarrow \frac{\ln(\sigma_{d}(t))}{\ln(t)} \quad\text{ as }\quad n \rightarrow \infty
\] 
as required. \QED}
\end{thm}

Figure \ref{fig:CircleDimensionComparisons} in Section \ref{sec:Applications} illustrates this result, showing an example of finite subspaces of a manifold $M$ approximating the quantity $\mathcal{F}\sigma_d(t)$ for its finite subsets.

We will now show the equivalent result for instantaneous spread dimension, which requires some intermediate results.
The following lemma is a straightforward application of Lebesgue's Bounded Convergence Theorem, see for example \cite[Theorem 3.26]{Axler2020:MeasureTheory}.

\begin{lem}\label{lem:integral_derivative_interchange}
Let $X$ be a compact space and let $f: \Rplus \times X \to \R$ be a continuous function. If $X$ is equipped with a finite measure $\mu$, and it for each $t \in \Rplus$ the partial derivatives $
\frac{\partial}{\partial t}f(t,x)$
exists, then
\[
\frac{\dee}{\dee t} \int\limits_{x \in X}f(t,x)\dee\mu(x) = \int\limits_{x \in X} \frac{\partial}{\partial t}f(t,x) \dee\mu(x)
\]
\end{lem}

\begin{lem}\label{lem:LemSigmaDifferentiable}
For $X$ a compact metric space with finite measure $\mu$ the function $\sigma_d(t):\Rplus \to \Rplus$ is differentiable and
\begin{equation}\label{eq:DerivativeOfSigma}
\frac{\dee}{\dee t}\sigma_d(t) = \int\limits_{x \in X} \frac{\int\limits_{y \in X} d(x,y)e^{-td(x,y)}\dee\mu(y) }{\big( \int\limits_{z \in X} e^{-td(x,z)}\dee\mu(z)  \big)^2} \dee\mu(x)
\end{equation}
\proof{Let $\varphi(t,x) = \int\limits_{y \in X}e^{-td(x,y)}\dee\mu(x)$, and let $\theta(t,x) = \frac{1}{\varphi(t,x)}$.

First we show that $\varphi(t,x)$ is continuous in $t$ and $x$, that is for any $t_n \rightarrow t$ as $n\rightarrow \infty$ and $x_m \rightarrow x$ 
as $m \rightarrow \infty$ we have $\lim\limits_{n,m \rightarrow \infty}\varphi(t_n, x_m) = \varphi(t,x)$.

Let $g_n(x) = \varphi(t_n,x)$. Note that for any monotone sequece $t_n\leq t_{n+1}$ the sequence $g_n(x)$ is also monotone. Since $X$ is compact, Dini's Theorem \cite[Theorem 7.13]{Rudin1976:PrinciplesOfMathematicalAnalysis} asserts that the convergence $g_n\rightarrow g$ is uniform. Hence by Lemma \ref{lem:ApostolLemma} we have 
\begin{align*}
\lim\limits_{n,m \rightarrow \infty}\varphi(t_n, x_m) =  \lim\limits_{n \rightarrow \infty} \lim\limits_{m \rightarrow \infty}\varphi(t_n, x_m) = \varphi(t,x)
\end{align*}
and hence $\varphi$ is continuous.

Next, we claim that $\theta(t,x)$ as also continuous in $t$ and $x$. Let $(a,b)$ be an open neighbourhood containing $t$, and let $r = \diam(M)$. Let $L = e^{-br}$, then note that that for all $t \in (a,b)$ and for all $x,y\in X$ we have $L \leq e^{-td(x,y)}$ and therefore we have
\begin{align*}
|\theta(t_n,x_m) - \theta(t,x)| = \Big| \frac{\varphi(t,x) - \varphi(t_n,x_m)}{\varphi(t_n,x_m)\varphi(t,x)} \Big| 
\leq \Big| \frac{\varphi(t,x) - \varphi(t_n,x_m)}{L^2} \Big| 
\end{align*}
and hence continuity of $\theta$ follows from continuity of $\varphi$.

Since the partial derivatives $\frac{\partial}{\partial t} e^{-td(x,y)}$ exist for all $t\in\Rplus$, it follows from Lemma \ref{lem:integral_derivative_interchange} that
\begin{align*}
\frac{\partial}{\partial t}\varphi(t,x) &= \int\limits_{y \in X} \frac{\partial}{\partial t} e^{-td(x,y)} = -\int\limits_{y \in X} d(x,y)e^{-td(x,y)} \
\end{align*}
By the reciprocal rule for differentiable functions, the partial derivatives $\frac{\partial}{\partial t}\theta(t,x)$ also exist with
\begin{align*}
 \frac{\partial}{\partial t} \theta(x,t)= \frac{-\frac{\partial}{\partial t}\varphi(t,x)}{\varphi(t,x)^2}= \frac{\int\limits_{y \in X} d(x,y)e^{-td(x,y)}}{\big(\int\limits_{y \in X} e^{-td(x,y)}\big)^2} 
\end{align*}
hence the result follows immediately from Lemma \ref{lem:integral_derivative_interchange}.
 \QED}
\end{lem}

\begin{thm}\label{thm:FiniteApproximation2}

Let $(M,d)$ be a compact metric space with a sequence of partitions $P_n$ satisfying the finite sampling condition. For each $n\in \N$ let $X_n \subset M$ be a set of representative points consisting of one point $x^n_i \in I_i^n$ for each subset in the partition $I_i^n \in P_n$. For each $t \in\Rplus$ we have
\[
\mathcal{G}\sigma_{d_n}(t) \rightarrow \mathcal{G}\sigma_d(t) \quad\text{ as }\quad n \rightarrow \infty
\] 
where $\sigma_{d_n}(t)$ denotes the spread of $X_n$.

\proof{By Lemma \ref{lem:SpreadApproximation}, it is enough to show that $\frac{\dee}{\dee t} \sigma_{d_n}(t) \rightarrow \frac{\dee}{\dee t}\sigma_M(t)$. It follows from Lemma \ref{lem:LemSigmaDifferentiable} that
\[
\frac{\dee}{\dee t}\sigma_{d_n} = \sum\limits_{x \in X_n} \frac{\sum\limits_{y \in X_n} d(x,y)e^{-td(x,y)}}{\Big(\sum\limits_{z \in X_n}e^{-td(x,z)}\Big)^2}
\]

Let $h_{n,m}(x)$ be defined
\[
h_{n,m}(x) = \frac{\sum\limits_{y \in X_n} d(x,y)e^{-td(x,y)}\frac{1}{n}}{\Big(\sum\limits_{z \in X_m}e^{-td(x,z)}\Big)^2\frac{1}{m^2}}
\]
and let $k_n(x) = \lim\limits_{m \rightarrow \infty} h_{n,m}(x)$. We need to show that for each $x$, the convergence $h_{n,m}(x) \rightarrow k_n(x)$ is uniform in $n$, that is, for each $\varepsilon>0$ there exists $M\in \N$ such that for all $m> M$ we have $|k_n(x) - h_{n,m}(x)|$ for all $n\in \N$.

Note that $\sum\limits_{y \in X_n} d(x,y)e^{-td(x,y)}\frac{1}{n} \rightarrow \int\limits_{y \in X} d(x,y)e^{-td(x,y)}\dee\mu(y)$ converges uniformly in $x$, using same argument as for $\varphi_n \to \varphi$ in the proof of Lemma \ref{lem:SpreadApproximation}. Letting $L = e^{-t r}$ where $r = \diam(M)$, then for any $\varepsilon>0$ we can pick $N\in \N$ such that for all $n>N$ 
\[
\Big|\sum\limits_{y \in X_n} d(x,y)e^{-td(x,y)}\frac{1}{n} - \int\limits_{y \in X} d(x,y)e^{-td(x,y)}\dee\mu(y)\Big|<\varepsilon L^2\]
Since for all $m\in \N$ we have $\sum\limits_{z \in X_n}e^{-td(x,z)}\frac{1}{m} \geq L$, it follows that for all $\varepsilon>0$ we have

\begin{align*}
\Big|\frac{\sum\limits_{y \in X_n} d(x,y)e^{-td(x,y)}\frac{1}{n}}{\Big(\sum\limits_{z \in X_n}e^{-td(x,z)}\frac{1}{m}\Big)^2} - \frac{\int\limits_{y \in X} d(x,y)e^{-td(x,y)}\dee\mu(y)}{\Big(\sum\limits_{z \in X_n}e^{-td(x,z)}\frac{1}{m}\Big)^2}\Big| &<\frac{\varepsilon L^2}{\Big(\sum\limits_{z \in X_n}e^{-td(x,z)}\frac{1}{m}\Big)^2}\\
&\leq \frac{\varepsilon L^2}{L^2} =\varepsilon
\end{align*}
and hence the convergence is uniform in $m$, and therefore by Lemma \ref{lem:ApostolLemma} the double limit exists, that is, for each $x \in X$
\begin{equation}\label{eq:eq_convergence}
\frac{\sum\limits_{y \in X_n} d(x,y)e^{-td(x,y)}}{\Big(\sum\limits_{y \in X_n}e^{-td(x,y)}\Big)^2\frac{1}{n}} \rightarrow \frac{\int\limits_{y \in X} d(x,y)e^{-td(x,y)}}{\Big(\int\limits_{y \in X} e^{-td(x,y)}\Big)^2} \quad\text{ as }\quad n\rightarrow\infty
\end{equation}
The expression (\ref{eq:eq_convergence}) can be see as a product of uniformly convergent sequences of functions $f_n(x) \rightarrow f(x)$, $g_n(x)\rightarrow g(x)$, where both sequences are uniformly bounded, and hence their products converge $f_n(x)g_n(x) \rightarrow f(x)g(x)$ uniformly in $x$. Therefore, the convergence (\ref{eq:eq_convergence}) is uniform in $x$.

Since this convergence is uniform, for any $\varepsilon>0$ we have $N\in \N$ such that for all $n>N$
\begin{equation}\label{eq_another_little_eq}
\Big|\frac{\sum\limits_{y \in X_n} d(x,y)e^{-td(x,y)}}{\Big(\sum\limits_{y \in X_n}e^{-td(x,y)}\Big)^2\frac{1}{n}} - \frac{\int\limits_{y \in X} d(x,y)e^{-td(x,y)}}{\Big(\int\limits_{y \in X} e^{-td(x,y)}\Big)^2}\Big|< \varepsilon
\end{equation}
For all $m\in \N$ we have
\begin{align*}&\Big|\sum\limits_{x \in X_m}\frac{\sum\limits_{y \in X_n} d(x,y)e^{-td(x,y)}}{\Big(\sum\limits_{y \in X_n}e^{-td(x,y)}\Big)^2\frac{1}{n}}\frac{1}{m} - \sum\limits_{x \in X_m}\frac{\int\limits_{y \in X} d(x,y)e^{-td(x,y)}}{\Big(\int\limits_{y \in X} e^{-td(x,y)}\Big)^2}\frac{1}{m}\Big| \\
 &\leq
\sum\limits_{x \in X_m}\Big|\frac{\sum\limits_{y \in X_n} d(x,y)e^{-td(x,y)}}{\Big(\sum\limits_{y \in X_n}e^{-td(x,y)}\Big)^2\frac{1}{n}} - \frac{\int\limits_{y \in X} d(x,y)e^{-td(x,y)}}{\Big(\int\limits_{y \in X} e^{-td(x,y)}\Big)^2}\Big| \frac{1}{m}\\
&< \sum\limits_{x \in X_m}\varepsilon \frac{1}{m} \quad\text{by (\ref{eq_another_little_eq})} \\
&=\varepsilon
\end{align*}
and hence this convergence is uniform in $m$, and therefore by Lemma \ref{lem:ApostolLemma} the joint limit exists and so
\[
\sum\limits_{x \in X_n}\frac{\sum\limits_{y \in X_n} d(x,y)e^{-td(x,y)}}{\Big(\sum\limits_{y \in X_n}e^{-td(x,y)}\Big)^2} \rightarrow \int\limits_{x \in X}\frac{\int\limits_{y \in X} d(x,y)e^{-td(x,y)}}{\Big(\int\limits_{y \in X} e^{-td(x,y)}\Big)^2}\dee\mu(x)
\]
as $n\rightarrow\infty$, hence for each $t \in \Rplus$ we have $\frac{\dee}{\dee t}\sigma_{d_n}(t) \rightarrow \frac{\dee}{\dee t} \sigma_X(t)$, as $n \rightarrow \infty$, as required.
\QED
}
\end{thm}
Figure \ref{fig:CircleDimensionComparisons} in Section \ref{sec:Applications} illustrates this result, showing an example of finite subspaces of a manifold $M$ approximating the instantaneous spread dimension of $M$.

\subsection{Estimating the Spread Dimension of Manifolds}

In order to apply the approximation results Theorem \ref{thm:FiniteApproximation1} and Theorem \ref{thm:FiniteApproximation2} to the case of Riemannian manifolds, we will now show that manifolds satisfy the finite samping condition.

The following lemma relates the Lebesgue integral to the integration of densities on manifolds, see, for example \cite[Chap. 16]{Lee2013:SmoothManifolds}. The result is well-known in the literature, but for lack of a definitive reference we sketch the proof.

\begin{lem}\label{lem:density_gives_radon}
For $M$ a Riemannian manifold, the Riemannian density $\omega$ yields a Radon measure $\mu$ such that for every Borel set $U\subset M$
\[
\mu(U) = \int_{U}\omega
\]
where the integration on the right hand side is integration of a density on a manifold.
\proof{The Riemannian density defines integration on the manifold, which yields a linear functional $\widetilde{\omega}:C(M) \rightarrow \R$ defined
\[
\widetilde{\omega}(f) = \int_{M} f \omega
\]
where linearity of $\widetilde{\omega}$ is demonstrated in, for example, \cite[Proposition 16.42]{Lee2013:SmoothManifolds}.

By the Riesz representation theorem -- see, for example \cite[Theorem 2.14]{Rudin1987:RealAndComplexAnalysis} -- the linear functional $\widetilde{\omega}$ uniquely corresponds with a Radon measure $\mu$ on $M$ satisfying
\[
\widetilde{\omega} (f) = \int_M f \dee\mu
\]
and hence for each $U\subset M$ we have
\[
\mu(U) = \int_U 1 \dee \mu 
=\int_U \omega 
\]
as required. \QED}
\end{lem}

The following result was stated by Kloeckner in answer to a question on Stack Overflow \cite{Kloeckner2016:StackOverflowAnswer} -- for the sake of completeness we give a detailed account of Kloeckner's proof.

\begin{lem}[Kloeckner]\label{lem:KloecknerLemma}
Let $M$ be an $n$-dimensional Riemannian manifold equipped with a Radon measure $\mu$ as defined in Lemma \ref{lem:density_gives_radon}. Any triangulation $\alpha: \mathcal{T} \rightarrow M$ is diffeomorphic to one satisfying $\mu( T_i) = \mu( T_j)$ for all $n$-simplices $T_i, T_j \in \mathcal{T}$.

\proof{ It follows from a result due to Moser \cite{Moser1965:VolumeElements} that if $\tau$ and $\theta$ are densities on $M$  such that 
\[
\int_M \tau = \int_M \theta
\]
then there exists a diffeomorphsm $\phi: M \to M$ such that for all $S \subset M$, we have
\[
\int_S \tau = \int_{\phi(S)}\theta
\]

Let $\alpha: \mathcal{T} \to M$ be a triangulation of $M$ such that the number of $n$-dimensional simplices in $\mathcal{T}$ is $k$.

Now let $f :M \to \R$ be a smooth function satisfying the property that for each $n$-dimensional simplex $T \in \mathcal{T}$ we have
\[
\int_{T} f\omega = \frac{\vol(M)}{k}
\]
To see that such a function $f$ exists, consider a constant function on $M$, and fixing the values on the boundary of each $n$-simplex $T \in \mathcal{T}$, allow the value $f$ to vary smoothly on the interior of each $T \in \mathcal{T}$ to give the required volume.

Now since we have 
\[
\int_M f\omega = \int_M \omega
\]
by the result of Moser, there exists a diffeomorphism $\phi: M \to M$ such that such that for each $T \in \mathcal{T}$ we have
\begin{align*}
\int_{\phi(T)} \omega & = \int_T f \omega  = \frac{\vol(M)}{k} 
\end{align*}
and hence it folllows from Lemma \ref{lem:density_gives_radon} that $\phi\circ \alpha: \mathcal{T} \rightarrow M$ is a triangulation of $M$ satisfying $\mu(T) =\frac{\vol(M)}{k}$ each $n$-simplex $T \in \mathcal{T}$, as required. \QED}
\end{lem}

\begin{thm}\label{thm:ManifoldsFiniteSampling}
Every compact Riemannian manifold satisfies the finite sampling condition, that is, admits a sequence of partitions satisfying the conditions of Definition \ref{def:FiniteSamplingCondition}.

\proof{ Let $M$ be an $n$-dimensional Riemannian manifold. We need to show that there exists a sequence of partitons of $M$ satisfying the criteria of Definition \ref{def:FiniteSamplingCondition}. 

First we note that every smooth manifold admits a triangulation, see for example \cite[p. 124]{Whitney1957:GeometricIntegration}.

Starting with any triangulation of $M$, barycentric subdivision of the $n$-simplices can be performed an aribtrary number of times.
After performing barycentric subdivision at every $n$-simplex, any pair of $n$-simplices sharing an $(n-1)$-simplex as a face can be recombined into one $n$-simplex, and therefore, by this process we can obtain a triangulation with any desired number of $n$-simplices. By Lemma \ref{lem:KloecknerLemma} any triangulation can be transformed into one such that all $n$-simplices have the same volume.

Furthermore, repeated barycentric subdivision will result in $n$-simplices with arbitrarily small diameter, see, for example \cite[p. 120]{Hatcher2002:AT}. Hence the process of subdividing triangulations and rescaling we can achieve a sequence of partitions satisfying the finite sampling condition.
\QED}
\end{thm}

\section{Examples with Real-World Data}\label{sec:Applications}

We will now give computations of spread and spread dimension for two real-world datasets, presented as a method of estimating the intrinsic dimension of those datasets. These results serve as proof of concept that the spread-based methods have real-world applications for intrinsic dimension estimation. All computations have been implemented in Python\footnote{Source code can be found at: \url{https://github.com/dk-gh/metric_space_spread}}, and in Appendix \ref{sec:appendix} we prove the correctness of efficient vectorised algorithms for computing spread and spread dimension.

Before considering examples, we briefly remark on some aspects of spread dimension as a method for intrinsic dimension estimation:

\begin{itemize}
\item Intrinsic dimension estimation techniques are broadly characterised as either global or local, where global techniques consider the whole dataset to infer the dimension, whereas a local technique is applied to small neighbourhoods only. We will see spread dimension is effective in measuring intrinsic dimension both globally and locally. This is noteworthy given the view that ``[a]ll the recent methods have abandoned the global approach since it is now clear that analyzing a dataset at its biggest scale cannot produce reliable results'' \cite[p. 3]{CampadelliEtAl2015:IntrinsicDimension}.

\item Spread dimension does not require data to be of the form $X \subset \R^n$. Although this type of data is the focus of the present work, there are many real-world non-Euclidean datasets where intrinsic dimension is a useful concept, for example, the fractal geometry of networks \cite{SongEtAl2007:FractalDimension}.

\item Computing spread dimension is \emph{manifold adaptive} \cite{Farahmand2007:ManifoldAdaptive}, which means the dimension of the ambient space does not affect the computational complexity of the spread dimension algorithms. The input of the algorithm is the pairwise distances of the points -- the complexity is therefore determined only by the number of points, and not the dimension. We will discuss practical significance of this in Remark \ref{rem:manifold_adaptive}. 

\item Spread dimension is \emph{robust against multiscaling} \cite{Camastra2016:IntrinsicDimension}. Many intrinsic dimension estimation techniques can give different results for data which differs only by a constant scale factor. The spread dimension approach is based on measuring how the spread grows across different scales, and is therefore inherently robust against multiscaling -- however, see Remark \ref{rem:robust_against_multiscaling} for a more detailed discussion.

\end{itemize}

Before considering real-world examples, we demonstrate the method using synthetic data, namely, points uniformly sampled from the unit circle in $\R^2$. Due to the idealised nature of this example, it gives a clear results that are easy to interpret, and we can compare the spread dimension of the finite samples to the exact formula of the spread dimension of the circle, which is known due to the following result of Willerton \cite[Theorem 7]{Willerton2013:MagnitudeOfSpheres}.

\begin{thm}[Willerton]\label{thm:SnExactSpread}
Let $(S^n, d)$ be an $n$-dimensional sphere of radius $1$, with intrinsic distance function $d$, where the distance between $x$ and $y$ is taken to be the length of the shortest path connecting $x$ and $y$.
Then for each $n\geq 1$, the spread $\sigma_d(t)$ is defined
\[
\sigma_d(t) = \begin{cases}
\frac{\pi t}{1-e^{-t \pi}}   & \text{ if } n=1 \\ \frac{2}{1+e^{-t \pi}}\prod\limits_{i=1}^{n/2}\big( (\frac{t}{2i-1})^2 +1 \big) & \text{ if } n \text{ even, and } n\geq 2 \\
\frac{\pi t}{1-e^{-t \pi}}\prod\limits_{i=1}^{(n-1)/2}\big( (\frac{t}{2i})^2 +1 \big) & \text{ if } n \text{ odd, and } n > 2 
\end{cases}
\]
\end{thm}

Willerton proves Theorem \ref{thm:SnExactSpread} for magnitude, but for homogenous metric spaces magnitude coincides with the spread \cite[Theorem 1]{Willerton2012:AsymptoticMagnitude}.

Using the exact formula for the spread of $(S^1,d)$ from Theorem \ref{thm:SnExactSpread}, we can derive the following exact formulae
\begin{equation}\label{eq:exact_formulae_circle}
\mathcal{G}\sigma_d(t) = 1 - \frac{\pi t}{e^{\pi t} - 1}\qquad \text{and} \qquad \mathcal{F}\sigma_d(t)  = 1 + \frac{\ln(\pi) - \ln(1-e^{-t\pi})}{\ln(t)}
\end{equation}

\begin{figure}[H]
\centering
    \advance\leftskip-4cm
        \advance\rightskip-4cm
\includegraphics[width=0.65\textwidth]{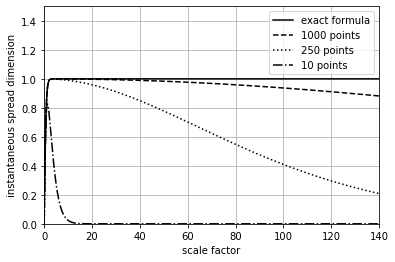}
\includegraphics[width=0.65\textwidth]{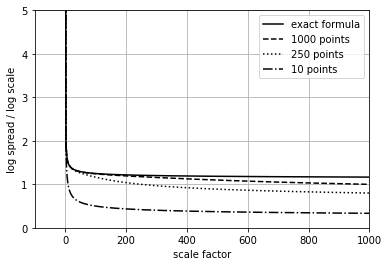}
\caption{\emph{(left)} The instantaneous spread dimension $\mathcal{G}\sigma_{d_n}(t)$ of various finite subsets of the circle plotted against the scale factor $t$. The value of $\mathcal{G}\sigma_d(t)$ obtained from the exact formula (\ref{eq:exact_formulae_circle}) is also plotted for comparison, illustrating how the instantaneous spread dimension is approximated by finite subsets.
\emph{(right)} The quantity $\mathcal{F}\sigma_{d_n}(t)$ for various finite subsets of the circle is computed for different numbers of points plotted against the scale factor $t$ for $t>1$. The value of $\mathcal{F}\sigma_d(t)$ obtained from the exact formula (\ref{eq:exact_formulae_circle}) is also plotted for comparison.} \label{fig:CircleDimensionComparisons}
\end{figure}

Figure \ref{fig:CircleDimensionComparisons} shows a comparison between the expressions for $\mathcal{G}\sigma_d(t)$ and $\mathcal{F}\sigma_d(t)$ for the circle against the corresponding values caluculated for finite subspaces with different numbers of points.

The dimension of $S^1$ can be inferred from either graph in Figure \ref{fig:CircleDimensionComparisons} -- using either $\mathcal{G}\sigma_d(t)$ or $\mathcal{F}\sigma_d(t)$ -- however the characteristics of these graphs are quite different. The way we use these quantities in practice to estimate dimension are outlined in the following heuristics.

In general, for any finite metric space we have $\mathcal{G}\sigma_d(0)=0$ and $\mathcal{G}\sigma_d(t)\rightarrow 0$ as $t \rightarrow \infty$, and hence the graph of $\mathcal{G}\sigma_d(t)$ will always have an identifiable peak, and Figure \ref{fig:CircleDimensionComparisons} shows typical behaviour where for large enough samples the instantaneous spread dimension reaches a plateau around the dimension of the space being sampled.

\begin{heuristic}
We use $\mathcal{G}\sigma_d(t)$ to estimate intrinsic dimension by looking for the peak or a plateau in the value of $\mathcal{G}\sigma_d(t)$. A long plateau is a stronger indication of the true instrinsic dimension than a short peak.
\end{heuristic}

The quantity $\mathcal{F}\sigma_d(t)$ is only shown for values of $t>1$, since $\mathcal{F}\sigma_d(t)$ is undefined at $t = 1$ and is negative for $0<t<1$. The values of $\mathcal{F}\sigma_d(t)$ fall from infinity as $t$ increases from values close to $t =1$.

\begin{heuristic}
We use $\mathcal{F}\sigma_d(t)$ to estimate intrinsic dimension by looking for the knee in the curve of $\mathcal{F}\sigma_d(t)$, as the value of $\mathcal{F}\sigma_d(t)$ falls from infinity.
\end{heuristic}

While Theorem \ref{thm:FiniteApproximation1} shows that using the quantity $\mathcal{F}\sigma_d(t)$ to estimate the dimension of $M$ is theoretically sound, in practice we will see that $\mathcal{G}\sigma_d(t)$ gives better results. A systematic approach to identifying an elbow or knee in a curve may be possible \cite{SatopaaEtAl2011:knees}, but there is no generally accepted definition of the ``knee'' and it operates more as a useful heuristic. In contrast to a knee, the peak value of $\mathcal{G}\sigma_d(t)$ is always unambiguously defined. The spread dimension $\mathcal{G}\sigma_d(t)$ also has other nice properties, like being robust against multiscaling -- see Remark \ref{rem:robust_against_multiscaling}.

\subsection{The Big-Five Personality Markers Dataset}

In this section we consider a dataset\footnote{Available at: \url{
https://openpsychometrics.org/_rawdata/IPIP-FFM-data-8Nov2018.zip
}} consisting of around one million responses to a psychometrics questionaire constructed from Goldberg's Big-Five factor markers  \cite{Goldberg1992:BigFive}, used to study human personality within the field of psychology.

The questionaire consists of fifty statements like: ``I often forget to put things back in their proper place'', or ``I talk to a lot of different people at parties'' to which the respondent assigns a numerical value in the range 1-5, where 1 means `disagree', 3 `neutral' and 5 `agree'.

The data therefore admits a natural representation as a one-million point subset $X \subset \R^{50}$, where each of the fifty dimensions corresponds with an individual question statement.

The Big-Five factor model of human personality asserts that human personality is broadly characterised by five independent traits: openness to experience; conscientiousness; extraversion; agreeableness; and neuroticism \cite[Chapter 1]{MatthewsEtAl2003:PersonalityTraits}. The five parameters along which personality varies manifest in empirical data having an intrinsic dimension of five. The intrinsic dimension of such data is typically identified using factor analysis, a linear dimension reduction technique commonly used in psychology and psychometrics.

The fact that five broadly consistent factors emerge from analysis of data from a variety of sources forms the empirical basis of Big-Five factor model, but there are alternative models proposing as few as three, or as many as 30 factors determining personality \cite[Chapter 1]{MatthewsEtAl2003:PersonalityTraits}. The intrinsic dimension of such empirical data is closely related to identifying the number of underlying factors, which can have fundamental significance to the development of theory  \cite{FabrigarEtAl1999:EvaluatingEFA,ZwickVelicer1986:ComparisonofFiveRules}.

Geometrically, factor analysis applied to data of the form $X\subset \R^m$, identifies linear subspace, or hyperplane of dimension $n$ within $\R^m$ such that the $X$ lies on a normal distribution in around that subspace \cite[Chapter 14]{HastieEtAl2009:ElementsOfStatisticalLearning}.

Since the data does not lie on the hyperplane, but is normally distributed around it, we would expect the raw data to exhibit a higher intrinsic dimension than the hyperplane. We will therefore consider the effects that smoothing the data has on the spread dimension. We use a nearest neighbour approach as a simple kernel smoothing method \cite[Chapter 6]{HastieEtAl2009:ElementsOfStatisticalLearning}.

\begin{defn}
Let $X \subset \R^n$ be a finite subset. For each  $x \in X$, let $N_k(x) \subset X$ be the set of $k$-nearest neighbours to $x$. For a choice of $k \in \N$, we define the following function $\xi^k: X \rightarrow \R^n$
\[\xi^k (x) = \frac{1}{k}\sum\limits_{y \in N_k(x)} y\]
where addition is defined point-wise in $\R^n$. The image of $\xi^k(X)$ is referred to as the \emph{knn smoothing} of $X$.
\end{defn}

\begin{figure}[H]
    \centering
    \advance\leftskip-4cm
        \advance\rightskip-4cm
\includegraphics[width=0.5\textwidth]{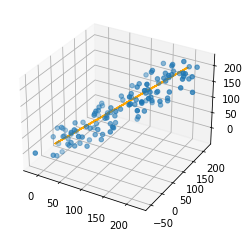}
\includegraphics[width=0.5\textwidth]{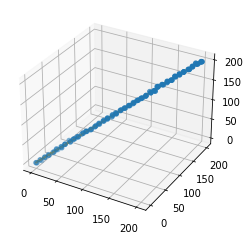}
\caption{The effect of smoothing on a set of points randomly distributed about some lower-dimensional hyperplane, in this case a line. Computing the spread dimension of the raw data \emph{(left)} would result in a dimension greater than one, whereas the smoothed data \emph{(right)} would have spread dimension tracking very closely to one.} 
\label{fig:KnnSmoothingDemo}
\end{figure}

Figure \ref{fig:KnnSmoothingDemo} shows the result of smoothing data that is normally distributed around a hyperplane.

\begin{figure}[H]
    \centering
    \advance\leftskip-4cm
        \advance\rightskip-4cm
\includegraphics[width=0.6\textwidth]{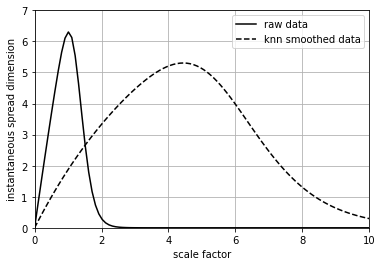}
\includegraphics[width=0.6\textwidth]{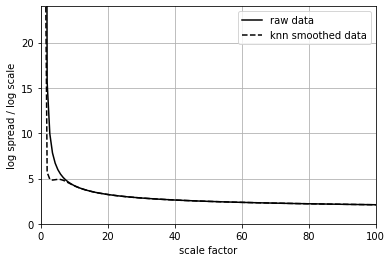}
\caption{For the closest 16,000 points around a randomly chosen point in the IPIP dataset, (\emph{left}) the instantaneous spread dimension $\mathcal{G}\sigma_d(t)$ plotted against scale factor $t$. (\emph{right}) The quantity $\mathcal{F}\sigma_d(t)$ plotted against scale factor $t$. Corresponding results for the knn-smoothed data where $k = 2400$ is plotted alongside the raw unsmoothed data.} \label{fig:ipip_local}
\end{figure}

With the Big-Five factor markers dataset we consider both global and local samples of the data. For the local samples, a point is chosen at random, along with its 16,000 nearest neighbours. We then generate a knn-smoothed data using a neighbouhbood size of 15\%, that is 2400 nearest neighbours. Figure \ref{fig:ipip_local} shows the results for a local sample.

\begin{figure}[H]
    \centering
    \advance\leftskip-4cm
        \advance\rightskip-4cm
\includegraphics[width=0.6\textwidth]{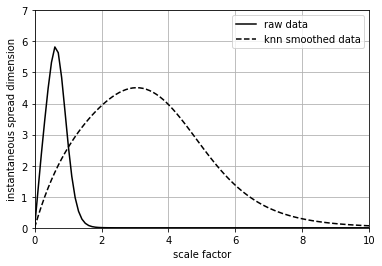}
\includegraphics[width=0.6\textwidth]{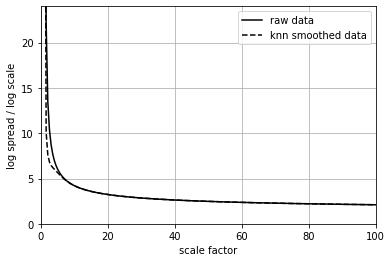}
\caption{For a random sample of 16,000 points from the IPIP dataset, (\emph{left}) the instantaneous spread dimension $\mathcal{G}\sigma_d(t)$ plotted against scale factor $t$. \emph{(right)} The quantity $\mathcal{F}\sigma_d(t)$ is plotted against scale factor $t$. Corresponding results for the knn-smoothed data where $k = 2400$ is plotted alongside the raw, unsmoothed data.}
\label{fig:ipip_global}
\end{figure}

We conduct the same analysis on global data, that is, a random selection of 16,000 points from the full dataset -- results shown in Figure \ref{fig:ipip_global}.

In both local and global sample, the smoothed data exhibits a maximum instantaneous spread dimension of five, rounded to the nearest integer. This is consistent with the intrinsic dimension inferred from factor analysis.

For the local data the spread dimension exhibits a clear elbow at five, but results for the global data are less easy to interpret.

\subsection{The Snake-Eyes Dice Image Dataset}

The \emph{Snake-Eyes} dataset\footnote{Created by Nicolau Werneck, available at: \url{https://github.com/nlw0/snake-eyes}} consists of one million computer-generated images of dice. This dataset was created in order to study translation and rotation in the context of image classification and manifold learning. Each image consists of $20 \times 20$ pixels, where each pixel is assigned a numerical greyscale value. Each image therefore admits a natural representation as a point in $\R^{400}$.

\begin{figure}[H]
\centering
\includegraphics[width=0.2\textwidth]{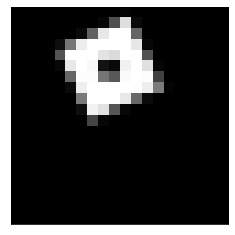}
\includegraphics[width=0.2\textwidth]{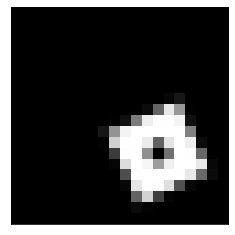}
\includegraphics[width=0.2\textwidth]{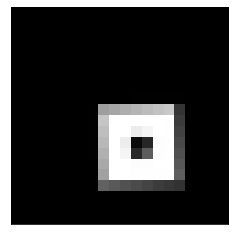}
\includegraphics[width=0.2\textwidth]{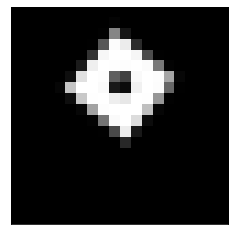}
\includegraphics[width=0.2\textwidth]{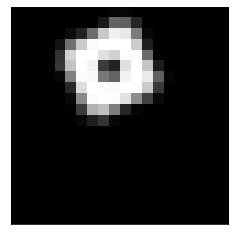}
\includegraphics[width=0.2\textwidth]{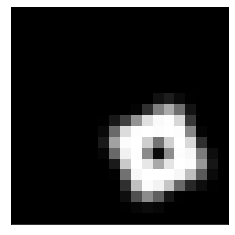}
\includegraphics[width=0.2\textwidth]{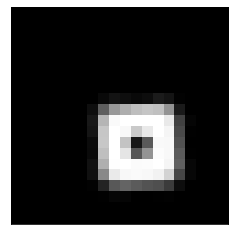}
\includegraphics[width=0.2\textwidth]{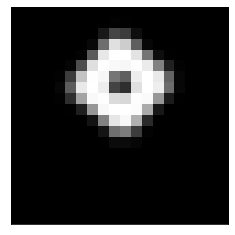}
\caption{\emph{(above)} Examples from the \emph{Snake Eyes dataset} represented as a $(20\times 20)$-pixel image. \emph{(below)} The corresponding images softened or smoothed using knn-smoothing.} \label{fig:snake_eyes}
\end{figure}

We consider the subset of images of a single dice showing the number one -- examples shown in Figure \ref{fig:snake_eyes}.

Locally, a dice on a surface has three degrees of freedom -- it can move in the $x$-axis, the $y$-axis or it can rotate about its centre. The underlying space $M\subset \R^{400}$ on which the data points lie is therefore isomorphic to $[0,1]^2 \times S^1$, and has intrinsic dimension $3$.

The Snake Eyes dataset is computer-generated and therefore it is not subject to random noise, however, the low-resolution of these images creates a roughness around the edges; a systematic error or artefact that we will see raises the intrinsic dimension. The idea that ``roughness of the boundary'' results in a measurable increase in dimension is a well-known heuristic in the applications of fractal geometry, for example in fractal analysis of medical imaging \cite{Kadi2008:Texture}. For this reason we will also consider image data subject to knn-smoothing, examples of which are also depicted in Figure \ref{fig:snake_eyes}.

\begin{rem}\label{rem:manifold_adaptive}
Instead of smoothing the data to eliminate the rough boundaries, we could simply use higher-resolution images. Note that using images with more pixels would not increase the time or space requirements on computing the spread and spread dimension, which are based only on the number of points -- that is, the method is \emph{manifold adaptive}. The significance of the manifold adaptive property is that the technique is not only applicable to small, low-resolution images.
\end{rem}

As with the previous section, we consider samples of points taken from the space globally, as well as locally, where a local sample consists of the $16,000$ points closest to a randomly chosen point.

\begin{figure}[H]
    \centering
    \advance\leftskip-4cm
        \advance\rightskip-4cm
\includegraphics[width=0.6\textwidth]{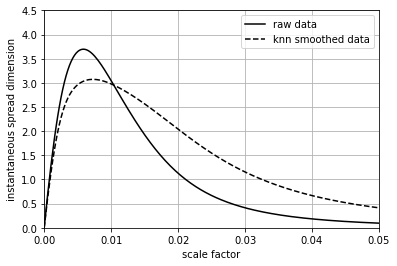}
\includegraphics[width=0.6\textwidth]{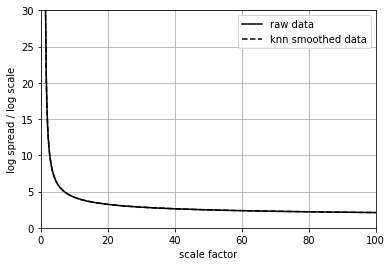}
\caption{For a random sample of 16,000 points from the Snake-Eyes dataset, (\emph{left}) the instantaneous spread dimension $\mathcal{G}\sigma_d(t)$ plotted against scale factor $t$. (\emph{right}) The quantity $\mathcal{F}\sigma_d(t)$ plotted against scale factor $t$.  Corresponding results for the knn-smoothed data is plotted alongside the raw unsmoothed data.}\label{fig:snake_results_global}
\end{figure}

\begin{figure}[H]
    \centering
    \advance\leftskip-4cm
        \advance\rightskip-4cm
\includegraphics[width=0.6\textwidth]{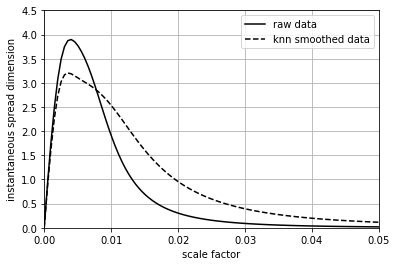}
\includegraphics[width=0.6\textwidth]{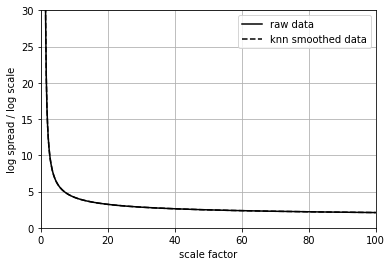}
\caption{For the closest 16,000 points around a randomly chosen point in the Snake-Eyes dataset, (\emph{left}) the instantaneous spread dimension $\mathcal{G}\sigma_d(t)$ plotted against scale factor $t$. \emph{(right)} The quantity $\mathcal{F}\sigma_d(t)$ plotted against scale factor $t$. Corresponding results for the knn-smoothed data is plotted alongside the raw unsmoothed data.}\label{fig:snake_results_local}
\end{figure}

In both the global case shown in Figure \ref{fig:snake_results_global}, and the local case shown in Figure \ref{fig:snake_results_local}, the instantaneous spread dimension of the smoothed data achieves a maximum value close to three. This is consistent with the intrinsic dimension of the underlying manifold. In both local and global cases the raw data exhibits a higher maximum instantaneous spread dimension of four. This may be due to the rough boundary of each dice image exhibiting an intrinsic dimension greater than one.

\begin{rem}\label{rem:robust_against_multiscaling}
It is clear from Figure \ref{fig:snake_results_global} and Figure \ref{fig:snake_results_local} that the quantity $\mathcal{F}\sigma_d(t)$ gives no clear results for the Snake Eyes data, and smoothing the data appears to have little discernible effect on the values of spread dimension. The reason from this is made clear by observing the scale at which meaningful values obtained: the scale at which spread captures meaningful geometric information is $t << 1$. We only consider the spread dimension for values of $t>1$, but at this scale the spread function is nearly constant with $\sigma_d(t) \approx 16,000$, and hence the graphs of spread dimension resemble the graph of $f(t) = \frac{16000}{\ln(t)}$, and contains little meaningful geometric information. 

This is an example of the scale skewing the accuracy of results; failure to be robust against multiscaling. On the other hand, the empirical results suggest that the instantaneous spread dimension gives results that are robust against multiscaling. This can be made a precise statement by observing that for any constant scale factor $A \in \Rplus$ we have
\[
\mathcal{G}\sigma_{Ad}(t) = \mathcal{G}\sigma_d(At)
\]
which follows from the chain rule. This implies that scaling the space merely has the effect of horizontally scaling the graph of the instantaneous spread dimension function, not changing its maximum value or shape.
\end{rem}

\subsection*{Acknowledgements}

I would like to thank Simon Willerton for useful comments and suggestions for this paper and the accompanying code. 

\bibliography{bibliography}

\begin{thebibliography}{10}

\bibitem{Kadi2008:Texture}
Omar~S. Al-Kadi and D.~Watson.
\newblock Texture analysis of aggressive and nonaggressive lung tumor {CE} {CT}
  images.
\newblock {\em IEEE Transactions on Biomedical Engineering}, 55(7):1822--1830,
  2008.

\bibitem{Apostol1974:Analysis}
Tom~M. Apostol.
\newblock {\em Mathematical Analysis}.
\newblock Addison-Wesley Publishing Company, second edition, 1974.

\bibitem{Axler2020:MeasureTheory}
Sheldon Axler.
\newblock {\em Measure, Integration and Real Analysis}.
\newblock Graduate Texts in Mathematics. Springer, 2020.

\bibitem{Camastra2016:IntrinsicDimension}
Francesco Camastra and Antonino Staiano.
\newblock Intrinsic dimension estimation: Advances and open problems.
\newblock {\em Information Sciences}, 328:26--41, 2016.

\bibitem{CampadelliEtAl2015:IntrinsicDimension}
P.~Campadelli, E.~Casiraghi, C.~Ceruti, and A.~Rozza.
\newblock Intrinsic dimension estimation: Relevant techniques and a benchmark
  framework.
\newblock {\em Mathematical Problems in Engineering}, 5:1--21, 2015.

\bibitem{Cunningham2015:LinearDimensionalityReduction}
John~P. Cunningham and Zoubin Gharamani.
\newblock Linear dimensionality reduction: Survey, insights and
  generalizations.
\newblock {\em Journal of Machine Learning Research}, 16:2859--2900, 2005.

\bibitem{Edgar2008:MeasureTopologyAndFractalGeometry}
Gerald Edgar.
\newblock {\em Measure, Topology, and Fractal Geometry}.
\newblock Undergraduate Texts in Mathematics. Springer, second edition, 2008.

\bibitem{FabrigarEtAl1999:EvaluatingEFA}
Leandre~R. Fabrigar, Duane~T. Wegener, Robert~C. MacCallum, and Erin~J.
  Strahan.
\newblock Evaluating the use of exploratory factor analysis in psychological
  research.
\newblock {\em Psychological Methods}, 4(3):272--299, 1999.

\bibitem{Farahmand2007:ManifoldAdaptive}
Amir~Massoud Farahmand, Csaba Szepesv\'{a}ri, and Jean-Yves Audibert.
\newblock Manifold-adaptive dimension estimation.
\newblock In {\em Proceedings of the 24th International Conference on Machine
  Learning ICML}, pages 265--272. ACM Press, 2007.

\bibitem{Goldberg1992:BigFive}
Lewis~R. Goldberg.
\newblock The development of markers for the {B}ig-{F}ive factor structure.
\newblock {\em Psychological Assessment}, 4(1):26--42, 1992.

\bibitem{GriepentrogEtAl2008:ABiLipschitzContinuous}
Jens~Andr\'{e} Griepentrog, Wolfgang H\"{o}ppner, Hans-Christoph Kaiser, and
  Joachim Rehberg.
\newblock A bi-{L}ipschitz continuous, volume preserving map from the unit ball
  onto a cube.
\newblock {\em Note di Matematica}, 1:177--193, 2008.

\bibitem{HastieEtAl2009:ElementsOfStatisticalLearning}
Trevor Hastie, Robert Tibshirani, and Jerome Friedman.
\newblock {\em The Elements of Statistical Learning: Data Mining, Inference,
  and Prediction}.
\newblock Springer Series in Statistics. Springer, second edition, 2009.

\bibitem{Hatcher2002:AT}
Allen Hatcher.
\newblock {\em Algebraic Topology}.
\newblock Cambridge University Press, 2001.

\bibitem{Higuchi1988:AnApproachToAnIrregularTimeSeries}
T.~Higuchi.
\newblock An approach to an irregular time series on the basis of the fractal
  theory.
\newblock {\em Physica D: Nonlinear Phenomena}, 31:277--283, 1988.

\bibitem{Kloeckner2016:StackOverflowAnswer}
Beno\^{i}t Kloeckner.
\newblock Answer to {S}tack {O}verflow question: Triangulation with simplices
  of same volume.
\newblock
  \url{https://mathoverflow.net/questions/237500/triangulation-with-simplices-of-same-volume},
  28/04/2016.

\bibitem{Knuth1976:BigO}
Donald~E. Knuth.
\newblock Big omicron and big omega and big theta.
\newblock {\em ACM Sigact News}, 8:18--24, 1976.

\bibitem{Lee1997:RiemannianManifolds}
John~M. Lee.
\newblock {\em Riemannian Manifolds: An Introduction to Curvature}.
\newblock Graduate Texts in Mathematics. Springer-Verlag New York Inc, 1997.

\bibitem{Lee2013:SmoothManifolds}
John~M. Lee.
\newblock {\em Introduction to Smooth Manifolds}.
\newblock Graduate Texts in Mathematics. Springer-Verlag New York Inc, second
  edition, 2013.

\bibitem{Leinster2013:Magnitude}
Tom Leinster.
\newblock The magnitude of metric spaces.
\newblock {\em Documenta Mathematica}, 18:857--905, 2013.

\bibitem{Leinster2021:Entropy}
Tom Leinster.
\newblock {\em Entropy and Diversity: An Axiomatic Approach}.
\newblock Cambridge University Press, 2021.

\bibitem{LiehrMassopust2020:OnTheMathematicalValidityOfHiguchi}
Lukas Liehr and Peter Massopust.
\newblock On the mathematical validity of the {H}iguchi method.
\newblock {\em Physica D: Nonlinear Phenomena}, 402, 2020.

\bibitem{MatthewsEtAl2003:PersonalityTraits}
Gerald Matthews, Ian~J. Deary, and Martha~C. Whiteman.
\newblock {\em Personality Traits}.
\newblock Cambridge University Press, second edition, 2003.

\bibitem{Meckes2015:Magnitude}
Mark~W. Meckes.
\newblock Magnitude, diversity, capacities, and dimensions of metric spaces.
\newblock {\em Potential Analysis}, 42:549--572, 2015.

\bibitem{Moser1965:VolumeElements}
J\"{u}rgen Moser.
\newblock On the volume elements on a manifold.
\newblock {\em Trans. Amer. Math. Soc.}, 120:286--294, 1965.

\bibitem{Rudin1976:PrinciplesOfMathematicalAnalysis}
Walter Rudin.
\newblock {\em Principles of Mathematical Analysis}.
\newblock McGraw-Hill, 3rd edition, 1976.

\bibitem{Rudin1987:RealAndComplexAnalysis}
Walter Rudin.
\newblock {\em Real and Complex Analysis}.
\newblock McGraw-Hill, 3rd edition, 1987.

\bibitem{SatopaaEtAl2011:knees}
Ville Satop\"{a}\"{a}, Jeannie Albrecht, David Irwin, and Barath Raghavan.
\newblock Finding a ``{K}needle'' in a haystack: Detecting knee points in
  system behavior.
\newblock In {\em 2011 31st International Conference on Distributed Computing
  Systems Workshops}, pages 166--171, 2011.

\bibitem{SongEtAl2007:FractalDimension}
Chaoming Song, Lazaros~K. Gallos, Shlomo Havlin, and Hern\'{a}n~A Makse.
\newblock How to calculate the fractal dimension of a complex network: the box
  covering algorithm.
\newblock {\em Journal of Statistical Mechanics: Theory and Experiment},
  2007(3), 2007.

\bibitem{Sutherland1975:IntroductionMetricTopology}
W.~A. Sutherland.
\newblock {\em Introduction to Metric and Topological Spaces}.
\newblock Oxford University Press, 1975.

\bibitem{TennenbaumEtAl2000:GlobalGeometricFramework}
Joshua~B. Tennenbaum, Vin de~Silva, and Jonh~C. Langford.
\newblock A global geometric framework for nonlinear dimensionality reduction.
\newblock {\em Science}, 290:2319--2323, 2000.

\bibitem{Maaten2009:DimensionalityReduction}
Laurens van~der Maaten, Eriv Postma, and Jaap van~den Herik.
\newblock Dimensionality reduction: A comparative review.
\newblock {\em Tilburg University Technical Report, TiCC-TR 2009-005}, 2009.

\bibitem{WaltEtAl2011:NumPy}
Stefan van~der Walt, S.~Chris Colbert, and Gael Varoquaux.
\newblock The {N}um{P}y array: A structure for efficient numerical computation.
\newblock {\em Computing in Science \& Engineering}, 13(2):22--30, 2008.

\bibitem{Whitney1957:GeometricIntegration}
Hassler Whitney.
\newblock {\em Geometric Integration Theory}.
\newblock Princeton University Press, 1957.

\bibitem{Willerton2013:MagnitudeOfSpheres}
Simon Willerton.
\newblock On the magnitude of spheres, surfaces and other homogeneous spaces.
\newblock {\em Geometriae Dedicata}, 168:291--310, 2014.

\bibitem{Willerton2013:Spread}
Simon Willerton.
\newblock Spread: A measure of the size of metric spaces.
\newblock {\em International Journal of Computational Geometry \&
  Applications}, 25(03):207--225, 2015.

\bibitem{Willerton2012:AsymptoticMagnitude}
Simon Willerton and Tom Leinster.
\newblock On the asymptotic magnitude of subsets of {E}uclidean space.
\newblock {\em Geometriae Dedicata}, 164:287--310, 2012.

\bibitem{Xiao1996:CartesianProduct}
Yimin Xiao.
\newblock Packing dimension, {H}ausdorff dimension and {C}artesian product
  sets.
\newblock {\em Mathematical Proceedings of the Cambridge Philosophical
  Society}, 120(3):535--546, 1996.

\bibitem{ZwickVelicer1986:ComparisonofFiveRules}
William~R. Zwick and Wayne~F. Velicer.
\newblock Comparison of five rules for determining the number of components to
  retain.
\newblock {\em Psychological Bulletin}, 99(3):432--442, 1986.

\end{thebibliography}

\appendix
\section{Vectorised Algorithms}\label{sec:appendix}

For $(X,d)$ a finite metric space with $n$ elments, we can fix an arbitrary order on the elements $x_1,x_2,...,x_n$, and define the distance matrix $D = (d_{i,j})$ to be the $n\times n$ matrix with entries $d_{i,j} = d(x_i,x_j)$.

We will now show that the spread and spread dimension of $(X,d)$ can be computed entirely in terms of matrix operations on $D$. This allows us to implement a vectorised algorithm in Python, where the entire computation is expressed as a sequence operations on NumPy arrays, which are highly optimised, enabling efficient computations \cite[p. 3]{WaltEtAl2011:NumPy}.

Assuming matrices $A = (a_{i,j})$ and $B = (b_{i,j})$ have the same dimensions, let $A\odot B$ denote their element-wise, or Hadamard product $A\odot B = (a_{i,j}\times b_{i,j})$. Element-wise squaring is denoted $A\odot A = A^{\odot 2}$, and the element-wise exponential $\exp(A) = (e^{a_{i,j}})$. Assuming $a_{i,j} \not= 0$ for all $i$ and $j$, then the element-wise inverse is defined $\inv(A) = \big(\frac{1}{a_{i,j}}\big)$.

We now show that the spread and its derivative can be expressed purely in terms of these operations performed on the distance matrix.

\begin{lem}\label{lem:vec1}
Let $(X,d)$ be a finite metric space with distance matrix $D$. If $J$ is the $n\times 1$-matrix with entries all 1, and $E_t = \exp(-tD)$ then
\[
\sigma_d(t) = J^{T} \cdot \big( \inv( E_t\cdot J )  \big)
\]

\proof{By straightforward calculation
\[
 E_t\cdot J= \left( {\begin{array}{cccc}
    e^{-t d_{1,1}} & e^{-td_{1,2}} & \cdots & e^{-td_{1,n}}\\
    e^{-td_{2,1}} & e^{-td_{2,2}} & \cdots & e^{-td_{2,n}}\\
    \vdots & \vdots & \ddots & \vdots\\
    e^{-td_{n,1}} & e^{-td_{m,2}} & \cdots & e^{-td_{n,n}}\\
  \end{array} } \right) \cdot \left( {\begin{array}{cccc}
    1\\
    1 \\
    \vdots \\
    1 
  \end{array} } \right) = \left( {\begin{array}{cccc}
     \sum\limits_{j= 1}^{n} e^{-t d_{1,j}}\\
    \sum\limits_{j= 1}^{n} e^{-t d_{2,j}} \\
    \vdots \\
    \sum\limits_{j= 1}^{n} e^{-t d_{n,j}} 
  \end{array} } \right) 
\]
and the sum of the element-wise inverses of the right hand side is exactly the expression for $\sigma_d(t)$, as required.
\QED}
\end{lem}

\begin{lem}\label{lem:vec2}
Let $(X,d)$ be a finite metric space with distance matrix $D$. If $J$ is the $n\times 1$-matrix with entries all 1, and $E_t = \exp(-tD)$ then
\[
\frac{\dee}{\dee t}\sigma_d(t) = J^{T} \cdot \Big( \big((D\odot E_t)\cdot J\big) \odot \big( \inv(E_t\cdot J)^{\odot 2}\big) \Big)
\]
\proof{ By straightforward calculation we have
\begin{align*}
D \odot E_t &= 
    \left( {\begin{array}{cccc}
    d_{1,1} & d_{1,2} & \cdots & d_{1,n}\\
    d_{2,1} & d_{2,2} & \cdots & d_{2,n}\\
    \vdots & \vdots & \ddots & \vdots\\
    d_{n,1} & d_{m,2} & \cdots & d_{n,n}\\
  \end{array} } \right)  \odot  \left( {\begin{array}{cccc}
    e^{-t d_{1,1}} & e^{-td_{1,2}} & \cdots & e^{-td_{1,n}}\\
    e^{-td_{2,1}} & e^{-td_{2,2}} & \cdots & e^{-td_{2,n}}\\
    \vdots & \vdots & \ddots & \vdots\\
    e^{-td_{n,1}} & e^{-td_{m,2}} & \cdots & e^{-td_{n,n}}\\
  \end{array} } \right) \\
  &=  \left( {\begin{array}{cccc}
    d_{1,1}e^{-t d_{1,1}} & d_{1,2}e^{-td_{1,2}} & \cdots & d_{1,n}e^{-td_{1,n}}\\
    d_{2,1}e^{-td_{2,1}} & d_{2,2}e^{-td_{2,2}} & \cdots & d_{2,n}e^{-td_{2,n}}\\
    \vdots & \vdots & \ddots & \vdots\\
    d_{n,1}e^{-td_{n,1}} & d_{m,2}e^{-td_{m,2}} & \cdots & d_{n,n}e^{-td_{n,n}}\\
  \end{array} } \right) 
\end{align*} 
and therefore
\begin{align*}
(D \odot E_t)\cdot J &= \left( {\begin{array}{cccc}
    d_{1,1}e^{-t d_{1,1}} & d_{1,2}e^{-td_{1,2}} & \cdots & d_{1,n}e^{-td_{1,n}}\\
    d_{2,1}e^{-td_{2,1}} & d_{2,2}e^{-td_{2,2}} & \cdots & d_{2,n}e^{-td_{2,n}}\\
    \vdots & \vdots & \ddots & \vdots\\
    d_{n,1}e^{-td_{n,1}} & d_{m,2}e^{-td_{m,2}} & \cdots & d_{n,n}e^{-td_{n,n}}\\
  \end{array} } \right) \cdot  \left( {\begin{array}{cccc}
    1\\
    1 \\
    \vdots \\
    1 
  \end{array} } \right)\\ 
  &= \left( {\begin{array}{cccc}
     \sum\limits_{j= 1}^{n} d_{1,j}e^{-t d_{1,j}}\\
    \sum\limits_{j= 1}^{n} d_{2,j}e^{-t d_{2,j}} \\
    \vdots \\
    \sum\limits_{j= 1}^{n} d_{n,j}e^{-t d_{n,j}} 
  \end{array} } \right) \\
\end{align*}
Also, we have
\[
\inv(E_t\cdot J)^{\odot 2} = \left( {\begin{array}{cccc}
     \frac{1}{\big(\sum\limits_{k= 1}^{n} e^{-t d_{1,k}}\big)^2}\\
    \frac{1}{\big(\sum\limits_{k= 1}^{n} e^{-t d_{2,k}}\big)^2}\\
    \vdots \\
    \frac{1}{\big(\sum\limits_{k= 1}^{n} e^{-t d_{n,k}}\big)^2}
  \end{array} } \right) 
\]
and therefore, taking the Hadamard product
\[
\Big((D \odot E_t)\cdot J\Big) \odot \Big(\inv(E_t\cdot J)^{\odot 2}\Big) = \left( {\begin{array}{cccc}
     \frac{\sum\limits_{j= 1}^{n} d_{1,j}e^{-t d_{1,j}}}{\big(\sum\limits_{k= 1}^{n} e^{-t d_{1,k}}\big)^2}\\
    \frac{\sum\limits_{j= 1}^{n} d_{2,j}e^{-t d_{2,j}}}{\big(\sum\limits_{k= 1}^{n} e^{-t d_{2,k}}\big)^2}\\
    \vdots \\
    \frac{\sum\limits_{j= 1}^{n} d_{n,j}e^{-t d_{n,j}}}{\big(\sum\limits_{k= 1}^{n} e^{-t d_{n,k}}\big)^2}
  \end{array} } \right) 
\]
and the sum of this column is precisely the formula for $\frac{\dee}{\dee t}\sigma_d(t)$ given in (\ref{eq:DerivativeOfSigma}), as required.
\QED}
\end{lem}

Lemma \ref{lem:vec1} and Lemma \ref{lem:vec2} demonstrate the correctness of the vectorised algorithm, which can be implemented in Python using NumPy array operations as follows.

\begin{center}
\begin{lstlisting}[style = Python]
import numpy as np
"""
D is a numpy array representing the distance matrix
for the metric space (X,d)
"""
def spread(D, t):
    """
    Returns the value of the scaled spread sigma_d(t)
    """
    J = np.ones(len(D))
    E = np.exp(-t*D)
    return np.sum(1/(E@J))

def instantaneous_spread_dimension(D, t):
    """
    Returns the value of the instantaneous growth
    dimension G(sigma_d(t))
    """
    J = np.ones(len(D))
    E = np.exp(-t*D)
    S = spread(D, t)
    return (t/S) * np.sum(((D*E)@J)/((E@J)**2))
\end{lstlisting}
\end{center}

\end{document}